\documentclass[a4paper]{article}

\author{Joaquim Ro\'e}

\usepackage[centertags]{amsmath}
\usepackage{amscd}
\usepackage{amsthm}
\usepackage{amssymb}
\usepackage{bbold}
\usepackage[dvips]{graphics}

\theoremstyle{definition}
\newtheorem{Def}{Definition}[section]
\theoremstyle{plain}
\newtheorem{Lem}[Def]{Lemma}
\newtheorem{Cor}[Def]{Corollary}
\newtheorem{Pro}[Def]{Proposition}
\newtheorem{Teo}[Def]{Theorem}

\theoremstyle{remark}
\newtheorem{Rem}[Def]{Remark}
\newtheorem{Exa}[Def]{Example}

\newcommand{\Ker}{\operatorname{Ker}}
\newcommand{\ord}{\operatorname{ord}}

\newcommand{\len}{\operatorname{length}}
\newcommand{\lcm}{\operatorname{l.c.m.}}

\newcommand{\edim}{\operatorname{edim}}

\newcommand{\lead}[1]{{#1}^{*}}
\newcommand{\val}[4]{{\widetilde{\operatorname{tr}}_{#3}^{#4}(#1|#2)}}
\newcommand{\Tr}[4]{{\operatorname{Tr}_{#3}^{#4}(#1|#2)}}
\newcommand{\Res}[4]{{\operatorname{Res}_{#3}^{#4}(#1|#2)}}
\newcommand{\TrZ}[4]{{\operatorname{SchTr}_{#3}^{#4}(#1|#2)}}
\newcommand{\ResZ}[4]{{\operatorname{SchRes}_{#3}^{#4}(#1|#2)}}
\newcommand{\HRes}[4]{{\operatorname{{\mathbb{Res}}}_{#3}^{#4}(#1|#2)}}
\newcommand{\tr}[4]{{\operatorname{tr}_{#3}^{#4}(#1|#2)}}
\newcommand{\res}[4]{{\operatorname{res}_{#3}^{#4}(#1|#2)}}
\newcommand{\hres}[4]{{\operatorname{{\mathbb{res}}}_{#3}^{#4}(#1|#2)}}
\newcommand{\Cl}{{\mathit{Cl}}}

\newcommand{\D}{{\bf D}}

\newcommand{\I}{{\cal I}}

\newcommand{\y}{{\mathfrak y}}

\newcommand{\M}{\mathfrak m}

\newcommand{\Z}{\mathbb{Z}}

\newcommand{\G}{\mathbb{G}}
\newcommand{\Q}{\mathbb{Q}}
\newcommand{\m}{{\bf m}}

\newcommand{\p}{{\bf p}}
\renewcommand{\y}{{\bf y}}
\renewcommand{\P}{\mathbb{P}}

\renewcommand{\O}{{\cal O}}
\renewcommand{\L}{{\cal L}}

\newcommand{\lmin}{\operatorname{\hat\ell_{min}}}
\newcommand{\lastr}{{r_{\max}}}

\newcounter{inleq}
\renewcommand{\theinleq}{\scshape \alph{inleq}}
\newcommand{\etag}[1]{\refstepcounter{inleq}\label{#1}\mbox{(\theinleq)}}
\newcommand{\eqrefprime}[1]{\mbox{(\ref{#1}$'$)}}

\begin{document}

\title{Limit linear systems and applications}
\maketitle

\begin{abstract}
A system of plane curves defined by prescribing $n$ points of multiplicity $e$ in general position is regular if $n\ge 4e^2$. 
The proof uses computation of limits of linear systems acquiring fixed divisors, an interesting problem in itself.
\end{abstract}

\section{Introduction}
\label{sec:results}

Linear systems defined by multiple points in the plane are a classical object of study, still not well understood.
Determining their regularity is one of the basic problems, equivalent to the solvability of bivariate Hermite interpolation problems and to the Riemann-Roch problem for rational surfaces. 
In spite of intense work devoted for decades to the question (\cite{Cil01}, \cite{Har02} and \cite{Mir99} are excellent overviews; even more recent results can be found in
\cite{BZ03}, \cite{CCMOF03}, 
\cite{Eva05}, \cite{Har03}, \cite{HR04}, \cite{Mig01}, \cite{Roe06},
\cite{Shu04}, \cite{Yan07}) it is far from settled and the main conjectures remain open.

Given $n$ points in the plane, and integers $e_1, \ldots e_n$, the curves of degree $d$ with multiplicity at least $e_i$ at the $i$-th point form a linear system $\mathcal L$ of dimension at least 
\begin{equation}
  \label{eq:vdim}
  \frac{d(d+3)}{2}-\sum_{i=1}^n \frac{e_i(e_i+1)}{2}.
\end{equation}
M. Nagata's famous conjecture of 1959 \cite{Nag59}, motivated by his solution to Hilbert's 14th problem, states that a nonempty linear system $\L$ defined by points in general position must have $d > \left(\sum e_i\right)/\sqrt{n}$ if $n>9$.
$\L$ is called \emph{regular} if it is empty or its dimension is given by \eqref{eq:vdim}.
In 1961, B. Segre \cite{Seg61} conjectured that a linear system $\L$ defined by points in general position is either regular or has a multiple curve in its base locus. The Harbourne\cite{Har86}--Hirschowitz\cite{Hir89} conjecture proposed in the 80's put further restrictions on the base curve;  C.~Ciliberto and R.~Miranda proved \cite{CM01} that it is equivalent to Segre's conjecture, and they both imply Nagata's conjecture.

In this work we focus on the equimultiple case $e_1=\cdots=e_n=e$, assuming that the points are in general position. The Harbourne--Hirschowitz conjecture implies in this case that $\L$ is regular if $n\ge 9$ (note that the conjecture is known to be true for $n\le 9$ \cite{Har85}, \cite{Har97}). 
Our main result is the following.

\begin{Teo}
\label{clm:regsmalle}
Let $e, n$ be positive integers with $n \ge 4e^2$. For general points $p_1, \ldots, p_n \in \P^2$, and for every $d$, the linear system $\mathcal L$ of curves of degree $d$ with multiplicity at least $e$ at each $p_i$ is regular.
\end{Teo}

For comparison purposes,  the only previous result which shows regularity for all $d$ when $n\ge f(e)$ for some function $f$ is due to J.~Alexander and A.~Hirschowitz \cite{AH00}, with $f \sim \exp(\exp(e))$.

Note also that regularity of $\L$ is known for all $d$ and all $n \ge 9$ if $e \le 42$ by recent work of M. Dumnicki \cite{Dum07}, so $\L$ is regular whenever $e \le \max\{42,\sqrt{n}/2\}$. 
Other known results for multiplicities small compared to the number of points support the weaker Nagata conjecture. Namely, in \cite{HR04} it is proved that an equimultiple system $\L$ with $d\le e\sqrt{n}$ is empty if $n\ge f(e)$, with $f(e)\sim e/2$ (L.~\'Evain \cite{Eva98} proved a similar result with $f(e)\sim 2e^2$).

It is also worth mentioning that regularity is known to hold for small nonequal multiplicities in some cases as well; to begin with, the aforementioned Alexander--Hirschowitz result holds for nonequimultiple systems (and even in higher dimension), and M.~Dumnicki--W.~Jarnicki \cite{DJ07} have proved regularity for all $d$ and all $n \ge 9$ if $e_i \le 11 \forall i.$ 
In a somewhat different spirit, S.~Yang \cite{Yan07} proved that, given an upper bound $e_i \le e \forall i$, there is a function $f(e)\sim e^2/\sqrt{6}$ such that, if $\L$ is regular for all $d$ and all $n \in [9, f(e)]$, then $\L$ is regular for all $d$ and all $n\ge 9$. 

\medskip

Let $k$ be an algebraically closed field of arbitrary characteristic,
and $X$ a smooth projective variety over $k$. Given an invertible
sheaf $L$ and a zero-dimensional scheme $Z\subset X$, a natural generalization of the preceding considerations is to ask about the regularity of the system $|L-Z|$  of effective divisors in $|L|=\P(H^0(X,\O_X(L)))$ containing $Z$. Such a system is regular if the natural linear map (restriction)
$$\Gamma(X,\O_X(L))\overset{\rho}{\longrightarrow} \Gamma(Z,\O_Z(L))$$
has maximal rank, as $|L-Z|=\P(\ker \rho)$.
It has revealed useful, when studying interpolation problems in
general position, to consider families of schemes $Z_t$ where the
position of the points supporting $Z_t$ varies with the parameter $t$.
Then one obtains a family of maps $\rho_t$ whose rank is lower semicontinuous
in $t$, so it is enough to find one value of the parameter, say $t=0$,
where the rank is 
maximal, to conclude that it is so for general $Z$. 
Several specialization techniques employed both classically (see
\cite{Pal03}, \cite{Ter16}, \cite{Nag59}) and recently (see
\cite{CM98a}, \cite{GLS98}, \cite{Har01}, \cite{Hir85},
\cite{Roe01a}, \cite{Roe01})  
rely on the fact that, if enough of $Z_0$
lies on a divisor $D$ of small degree, then all divisors in 
$|L-Z_0|$ must contain $D$, and subtracting $D$ gives a linear system of 
the same dimension with smaller degree and smaller $Z$;
then one hopes to show maximal rank inductively. A systematic
use of specialization
to divisors is sometimes called ``the Horace method'' after \cite{Hir85}.

The drawback to this method is that
if $\rho_0$ does not have maximal rank, it just gives a weak bound for
the actual behaviour in general position. Alexander-Hirschowitz
\cite{AH92a} (resumed and refined by the same authors in \cite{AH95} and
\cite{AH00}, by Mignon in \cite{Mig00a}, \cite{Mig00} and
\cite{Mig01}, by Chandler in \cite{Cha01} and by 
\'Evain in \cite{Eva??}, \cite{Eva07}) and
Ciliberto-Miranda in \cite{CM98a} (resumed and refined by the same authors in
\cite{CM00} and \cite{CM05} and by Buckley-Zompatori in \cite{BZ03})
have shown a way around this obstacle.  Denote by 
$$\edim |L-Z|=\max\{\dim \Gamma(X,\O_X(L))- \len Z-1, -1\}$$
the expected dimension. The idea is to 
consider the limit of $|L-Z_t|$ when $t$ tends to $0$ (in the
Grassmannian of $|L|$), and to construct a suitable
``intermediate'' interpolation problem
$$\Gamma(X,\O_X(L-mD))\overset{\rho_0'}{\longrightarrow} 
\Gamma(Z_0',\O_{Z_0'}(L-mD)),$$
in the sense that 
$$
\lim_{t\mapsto 0} |L-Z_t| \,\subset\, mD+|(L-mD)-Z_0'| \,\subset\, |L-Z_0|
$$
and $\edim |L-Z_t|=\edim |(L-mD)-Z_0'|$. Then
it is enough that $\rho_0'$ have maximal rank (rather than
$\rho_0$) to deduce that $\lim_{t\mapsto 0} |L-Z_t| = mD+|(L-mD)-Z_0'|$ 
and prove maximal rank for schemes in general position.

Here we explicit a method to systematically construct such
intermediate problems, in the particular case that only one point of
the support of $Z_t$ varies with $t$, and study its range of
applicability. 
The computation of limit linear
systems is interesting on its own, in addition to our original
motivation, for instance for the computation of limits in Hilbert
schemes \cite{Eva07}, and for adjacency of equisingularity types \cite{AR05}.
Hence we are interested in presenting the method in its natural general
setting; also, even for the applications to the plane, we actually need
to compute limit systems on other rational surfaces $X$.

Our approach is a generalization of the differential
Horace method as presented by \'Evain in \cite{Eva??}. 
In \cite{Eva07} \'Evain gives a further
generalization that allows several points to move, and even though his
statements deal with vertical translations of monomial schemes only,
our definition of the intermediate problem can be
implicitly found in his proofs. 
However, the methods of \cite{Eva??} and \cite{Eva07} don't directly apply
to general families (or even to families of monomial schemes moving
non-vertically) because
not every family allows an intermediate system of the same expected 
dimension as the original. It
may even happen that the limit system is not
determined by the condition of containing a subscheme.
We have identified the obstructions to the existence 
of such intermediate systems for a general family of zero-dimensional schemes
with one moving point, as elements in certain ideal quotients 
(proposition \ref{clm:correct}). 
Such obstructions did not appear in the method of \cite{Eva07},
because they in fact vanish for vertical translations of a monomial scheme.


To effectively apply the method to a particular family of zeroschemes,
some nontrivial algebraic computations are needed.
Here we restrict to families of monomial zeroschemes 
(moving non-vertically) and their
projections by blowing down. In a few cases we can then use the
computations of \cite{Eva07}. We expect however to exploit the
generality of the method in the future, as the
knowledge of obstructions should help in the search of useful
specializations.

\subsection{Limit linear systems}



Let $C$ be a quasi-projective smooth curve over $k$ and let $Z$ be a
subscheme of $X \times C$ which is flat and finite over $C$. 
The dimension of $|L-Z_t|$ is an upper semicontinuous function of $t\in C$ 
(with the Zariski topology of $C$). 
Thus there is an open set $U\subset C$ where $\dim
|L-Z_t|$ is minimal and constant, say $d$. This gives a morphism
to the Grassmannian of $d$-dimensional linear
subspaces of $|L|$,
\begin{align*}
  U&{\rightarrow} \G(|L|,d)\\
t&\mapsto |L-Z_t|
\end{align*}
which can be extended to the whole of $C$ because
the Grassmannian is projective. For $t_0 \not \in U$ we denote
$\lim_{t\mapsto t_0}|L-Z_t|$ the image of $t_0$ by the extension of
the morphism above to $C$. Let
$\mathcal I_t$ be the ideal sheaf of $Z_t$.
Assume that for some $t_0 \in C$ there exists a prime
divisor $D \subset X$ such that
$\rho_{t_0}^D:H^0(D,\O_D(L))\rightarrow H^0(D\cap Z_{t_0},\O_{D\cap Z_{t_0}}(L))$
is injective; then $D$ is a fixed part of $|L-Z_{t_0}|$. 
The residual linear system after subtracting $D$, which has the same
dimension as $|L-Z_{t_0}|$, is $|(L-D)-\tilde
Z_{t_0}|$, with $\tilde Z_{t_0}$ defined by the residual exact sequence
$$0\rightarrow \I_{\tilde Z_{t_0}}\rightarrow
\I_{Z_{t_0}}\rightarrow \I_{Z_{t_0}}\otimes \O_D \rightarrow 0.$$
In order to use the special member $Z_{t_0}$ to prove regularity 
of general $|L-Z_t|$, this residual should have the same
expected dimension as the original systems. But if $\rho_{t_0}^D$ is
not surjective, the expected dimension
will jump (in the language of Horace methods, the specialization is
not adjusted).

The jump in expected dimension comes from specializing ``too much'' of
$Z_{t_0}$ onto $D$. Now, the idea of \cite{AH00} is, roughly speaking,
to take $(t-t_0)^p=0$ for some $p>1$, so that $Z_t \cap D$ is big
enough to have $\rho_{t}^D$ injective, but not as much as
$Z_{t_0}\cap D$, so that $\rho_{t}^D$ can be adjusted.
Generalizing this idea as in \cite{Eva07}, if
$\p=(p_1, \dots, p_m) \in \Z^m$ is 
a non-increasing sequence, one takes more and more special positions given by
$(t-t_0)^{p_1}=0$, $(t-t_0)^{p_2}=0$, etc.   Then
one needs an adequate notion of residual to bound $\lim_{t\mapsto
  t_0}|L-Z_t|$ and show that these intuitions correspond rigorously to
an actual phenomenon. In section \ref{sec:method} 
iterated trace $\Tr{Z_t}{D}{\p}{i}$ and 
residual $\Res{Z_t}{D}{\p}{i}$ ideals are defined, providing
intermediate systems, and we show that
under suitable hypotheses, including (but not restricted
to) the specializations of monomial ideals of \cite{AH92a},
\cite{AH00}, \cite{Mig00a}, and \cite{Eva07},
they do have the same expected dimension as $|L-Z_t|$. 

Suppose that $Z=(Z_{\operatorname{fix}}\times C) \cup Y \subset X \times C$, 
where $Z_{\operatorname{fix}}\subset X$ is a 
fixed zero-dimensional scheme and $Y$ is irreducible, finite and
flat over $C$. In other words, $Z$ has a fixed and a moving part,
and the moving part is supported at a single (possibly moving) point of
$X$. Assume that $Y_{t_0}$ is supported at a point on a prime dvisor
$D$, and $\rho_{t}^D$ is injective. Denoting by
$\TrZ{Z_t}{D}{\p}{i}\subset D$ and $\ResZ{Z_t}{D}{\p}{i} \subset X$
the zero-dimensional subschemes defined by the trace and residual
ideals, the main result on limit linear systems is the following:


\begin{Teo}
\label{clm:global}
Let $Z=(Z_{\operatorname{fix}}\times C) \cup Y$ be as above, and take a sequence
$\p=(p_1, \dots, p_m) \in \Z^m$ with $p_1 \ge p_2 \ge \dots \ge p_m
\ge 1$. 
\begin{enumerate}
\item If for $1 \le i \le m$, the map
$$ \rho_{i}^D:H^0(D,\O_D(L-(i-1)D))\longrightarrow
H^0(\TrZ{Z_t}{D}{\p}{i},\O_{\TrZ{Z_t}{D}{\p}{i}}(L))$$
is injective, then
$\lim_{t\mapsto t_0}|L-Z_t| \subset mD+|(L-mD)-\ResZ{Z_t}{D}{\p}{m}|$.
\item Given $Z$ and $p_1, \dots, p_i$, $i<m$, there exists $1\le
  q_i$ such that, if for all $i=1, \dots m-1$, $ \rho_{i}^D$ is
  bijective, $p_{i+1}\le q_i$ and the restriction map
  $H^0(X,\O_X(L-(i-1)D))\rightarrow H^0(D,\O_D(L-(i-1)D))$ has maximal
  rank,  then 
  $\edim |L-Z_t|= \edim |(L-mD)-\Res{Z_t}{D}{\p}{m}|$.
\end{enumerate}
\end{Teo}

The results of section \ref{sec:method} are in fact slightly more
general, since we allow for singular and
reducible fixed divisors $D=D_1+ \dots +D_k$, and each component may
appear with a different multiplicity in $|L-Z_{t_0}|$. The first claim
of theorem \ref{clm:global} is a natural 
generalization of \'Evain \cite{Eva??}, more or less implicit in
\cite[Theorem 1]{Eva07}, but the 
second is to our knowledge entirely new, since the methods of
\cite{Eva07} (which give $q_i=p_i-1$ in the case of vertically
translated monomial schemes) do not apply in the general setting.

Whereas the bijectivity hypotheses in the second claim of 
theorem \ref{clm:global} are adjustment requirements (depending on the
global geometry of $D$), the hypotheses on $\p$ 
are of a new local kind: they force
that the obstructions mentioned above and specified in
proposition \ref{clm:correct} actually vanish for the given
specialization. 
The exact value of $q_i$ may be found in the proof of
corollary \ref{clm:decreasp}, but in the applications it will be
advantadgeous to apply results of section \ref{sec:homog}, where
sequences $\p$ are analyzed with respect to the valuative properties
of $Z$ and $D_i$. 

The proof of theorem \ref{clm:regsmalle} is based on theorem \ref{clm:global}. 
However, we don't actually compute the schemes $\TrZ{Z_t}{D}{\p}{i}\subset D$
and $\ResZ{Z_t}{D}{\p}{i} \subset X$; instead, we give bounds for them and 
use the second part of \ref{clm:global} to make sure that the expected
dimension is preserved and so the regularity of the limit system
proves regularity of the general ones.

Once theorem \ref{clm:regsmalle} is known, and using theorem \ref{clm:global},
the following result of \'Evain \cite{Eva07} (see also Ciliberto-Miranda 
\cite{CM06}) can be quickly proved:

\begin{Teo}\label{clm:squares}
  Assume that the characteristic of the base field $k$ is zero.
  Let $n=s^2$ be a square, and let $e\ge 1$ be an integer. 
  The scheme formed by $n$ distinct points in general position in
  $\P^2$ with multiplicity $e$ has maximal rank in all degrees.
\end{Teo}

An interesting feature of the proof of \ref{clm:squares} is that
it shows that
the method can still be useful when obstructions do appear, to 
prove emptyness of a linear system.
Also, it may be worth noting that it uses
the same specialization Nagata used in \cite{Nag59}
to prove that, for every square $n=s^2>9$ and every integer $e$ there are no
curves of degree $se$ with multiplicity at least
$e$ at $n$ general points.


\subsection*{Acknowledgements}

Many ideas underlying the present approach can be found in previous
articles by J.~Alexander-A.~Hirschowitz, L.~\'Evain,   and
T.~Mignon. Section \ref{sec:formal} was developed  
during my stay in Nice in 2001, and owes much to conversations with
A.~Hirschowitz, to whom I would like to express the most sincere gratitude.

I also want to thank E. Casas-Alvero for introducing me to the use of
clusters and Enriques diagrams, which have proved so useful in the
study of linear systems, and
M.~Alberich-Carrami\~nana, C.~Ciliberto, B.~Harbourne,
S.~Kleiman, R.~Miranda, and especially L.~\'Evain, for many
interesting conversations that 
greatly influenced this work.

\section{Algebraic approach to intermediate linear systems}
\label{sec:method}

Section \ref{sec:formal} contains the algebraic 
local analysis of the behaviour of a linear system moving in a
1-dimensional family that acquires a base divisor (possibly reducible
with multiple components) in a special position. 
In section \ref{sec:ncond} we determine under which conditions the
intermediate linear system given by our method coincides with the 
expected limit linear system, and we prove theorem \ref{clm:global}.
To effectively apply the results of this section, some computations
are needed which tend to be nontrivial. 
Section \ref{sec:sqcomp} shows one such computation, needed for the 
proof of theorem \ref{clm:squares}.
Part of these computations can be arranged in a systematic way, 
and we obtain sufficient conditions under which the intermediate 
linear system coincides with the expected limit system, for homogeneous
or monomial families, in sections \ref{sec:homog} and \ref{sec:monomial}.

\subsection{Higher order traces and residuals}
\label{sec:formal}

Let $R$ be an integral $k$-algebra, and consider $R_t=R \otimes
k[[t]]$. Given $f_t\in R_t$, denote $f_0\in R$ its image by the
obvious morphism $t \mapsto 0$. Similarly, for an ideal $I_t$ in
$R_t$, $I_0=(I_t+(t))/(t) \subset R_t/(t) \cong R$.

We define higher order traces and residuals of $I_t$ on divisors
$y=0$, in the spirit of \cite{Eva07}. Loosely speaking, if the $I_t$
define a family of schemes, we want to consider the trace on the
divisor $y=0$ of the special member given by $t^p=0$, and compute the
residual family (over $k[[t]]/(t^p)$) in a way that allows to
consider the trace on another 
divisor $z=0$ of every special residual given by $t^q=0$, $q\le p$, etc.

Given an ideal $I_t \subset R_t$, an element $y \in R$ and an integer
$p\ge 1$, consider the following ideals:
\begin{align*}
\Tr{I_t}{y}{p}{}= & \frac{\left((I_t+(y)):t^{p-1}\right)_0}{(y)}
\subset R/(y),\\ \HRes{I_t}{y}{p}{} =& (I_t+(t^{p})):y_1
\subset R_t,\\ \Res{I_t}{y}{p}{}= & \left((I_t+(t^{p})):y\right)_0
\subset R.
\end{align*}
Note that 
$\Res{I_t}{y}{p}{}=\left(\HRes{I_t}{y}{p}{}\right)_0$, and there are
inclusions $\Tr{I_t}{y}{1}{}\subset \Tr{I_t}{y}{2}{} 
\subset \dots$, and $\Res{I_t}{y}{1}{}\supset \Res{I_t}{y}{2}{} \supset
\dots \supset I_0$. 

More generally, given
sequences $\p=(p_1, \dots, p_m) \in \Z^m$ with $p_1 \ge p_2 \ge
\dots \ge p_m \ge 1$, $\y=(y_1, \dots, y_m) \in R^m$, denote
$\Tr{I_t}{\y}{\p}{1}=\Tr{I_t}{y_1}{p_1}{}$, 
$\Res{I_t}{\y}{\p}{1}=\Res{I_t}{y_1}{p_1}{}$ and 
$\HRes{I_t}{\y}{\p}{1}=\HRes{I_t}{y_1}{p_1}{}$; then for every integer
$1< i \le m$ define iteratively the following ideals:
\begin{align*}
\Tr{I_t}{\y}{\p}{i}= & \Tr{{\HRes{I_t}{\y}{\p}{i-1}}}
 {y_i}{p_i}{} \subset R/(y_i),\\
\HRes{I_t}{\y}{\p}{i}= & \HRes{{\HRes{I_t}{\y}{\p}{i-1}}}
 {y_i}{p_i}{}\subset R_t,\\
\Res{I_t}{\y}{\p}{i}= & \Res{{\HRes{I_t}{\y}{\p}{i-1}}}
 {y_i}{p_i}{} \subset R.
\end{align*}
If $y_i=y \, \forall i$ we  write 
$\Tr{I_t}{y}{\p}{i}=\Tr{I_t}{\y}{\p}{i}$, 
$\HRes{I_t}{y}{\p}{i}=\Tr{I_t}{\y}{\p}{i}$ and
$\Res{I_t}{y}{\p}{i}=\Tr{I_t}{\y}{\p}{i}$, and  
if $p_i=p \, \forall i$  we  write 
$\Tr{I_t}{\y}{p}{i}=\Tr{I_t}{\y}{\p}{i}$, etc. Note that
\begin{equation}
\label{eq:altres}
  \HRes{I_t}{\y}{\p}{i}=
\left(I_t+(t^{p_1}, y_1t^{p_{2}},\ \dots\ ,
y_1 y_2 \cdots y_{i-1}t^{p_{i}})\right): (y_1 y_2 \cdots y_{i}).
\end{equation}
Sometimes we shall also write $\HRes{I_t}{y}{\p}{0}=I_t$ and
$\Res{I_t}{y}{\p}{0}=I_0$. 


\medskip

Proposition \ref{clm:formal} below (or rather its immediate corollary
\ref{clm:formalah}) is a natural generalization of theorem 14 in
\cite{Eva??} (proved in \cite{Eva07} for products of ideals in
products of rings), which in turn refines proposition 8.1 of \cite{AH00}. It
justifies the definitions given
so far, and it will imply the first part of
theorem \ref{clm:global}. 
The reader may notice that the method of
proof is essentially the same used by \'Evain in \cite{Eva07}, theorem 1,
for the particular case that $R$ is a power series ring, one of
the variables is $y=y_i \, \forall i$, and $I_t$ is a monomial ideal of $R$
translated ``vertically'', i.e., by $y\mapsto y+t$. For this
particular case, equivalent definitions to the ones above can be found
in the proof of theorem 1 of \cite{Eva07} (in 
particular, $\HRes{I_t}{y}{\p}{i}$ is called
$J_{p_1:\dots:p_i:}$ there).

For a $k$-linear subspace $V \subset R$ and $y\in R$, let 
$\Res{V}{y}{}{}=\{v \in R \ | \ vy \in V\}$. 

\begin{Pro}
\label{clm:formal} Let $V \subset R$ be a $k$-linear subspace, and $I_t
\subset R_t$ an ideal such that $R_t/I_t$ is flat over 
  $k[[t]]$. Let $m\ge 0$ and let $\p=(p_1, \dots, p_m) \in \Z^m$ with
  $p_1 \ge p_2 \ge \dots \ge p_m \ge 1$, and $\y=(y_1, \dots, y_m) \in
  R^m$ be given sequences. Let $W=\{f \in R | \exists f_t\in  V
  \otimes k[[t]] \cap I_t \text{ with } f_0=f\}$. 
  If for $1 \le i \le m$, the canonical map
$$ \frac{\Res{V}{y_1 \cdots y_{i-1}}{}{}}{\Res{V}{y_1 \cdots
    y_{i-1}}{}{} \cap (y_i)} \longrightarrow
\frac{R/(y_i)}{\Tr{I_t}{\y}{\p}{i}}$$
is injective, then
$W \subset y_1 \dots y_m  \Res{I_t}{\y}{\p}{m}$.
\end{Pro}

\begin{proof}
Let $f_t \in V \otimes k[[t]] \cap I_t$. 
If 
\begin{equation}
  \label{eq:formres}
f_t\in (t^{p_1}, y_1t^{p_{2}},\ \dots\ ,
y_1 y_2 \cdots y_{m-1}t^{p_{m}},y_1 y_2 \cdots y_{m}),  
\end{equation}
i.e.,  if $f_t=g_t y_1 y_2 \cdots y_{m}+h_t$ for some $h_t \in
\linebreak[0] (t^{p_1}, y_1t^{p_{2}},\ \dots\ , y_1 y_2 \cdots
y_{m-1}t^{p_{m}})$, then (\ref{eq:altres}) implies 
$g_0 \in \Res{V}{y_1 y_2 \cdots y_{m}}{}{}\cap \Res{I_t}{\y}{\p}{m}$, 
and therefore $f_0\in y_1 \dots y_m  \Res{I_t}{\y}{\p}{m}$ as claimed. So it will
be enough to prove (\ref{eq:formres}).

We use induction on $m$. For $m=0$, one has trivially $f_t \in
(\prod_{i=1}^0 y_i) = R_t$. Assume now $m>0$ and $$f_t\in (t^{p_1}, y_1t^{p_{2}},\ \dots\ ,
y_1 y_2 \cdots y_{m-2}t^{p_{m-1}},y_1 y_2 \cdots y_{m-1}).$$
Denoting $p=p_m$,
and taking into account that $p_1 \ge p_2 \ge \dots \ge p_m \ge
1$, this implies $f_t \in (y_1 y_2 \cdots y_{m-1},t^{p})$, i.e., $f_t=g_t
y_1 y_2 \cdots y_{m-1}+h_t t^{p}$, where we may further assume that
$g_t=G_0+G_1t+\dots+G_{p-1}t^{p-1}$, with $G_j \in R$, $j=0,
\dots, p$. Denote $\bar G_j$ the class of $G_j$ in $R/(y_m)$; we
want to see that $\bar G_0= \dots= \bar G_{p-1}=0$.

The inclusions $\Tr{\HRes{I_t}{\y}{\p}{m-1}}{y_m}{0}{}\subset
\Tr{\HRes{I_t}{\y}{\p}{m-1}}{y_m}{1}{}  \subset \dots$
together with the hypothesis in the case $i=m$ tell us that,
for every $j=1, \dots, p$, the map $$
\varphi_j:\frac{\Res{V}{y_1 y_2 \cdots y_{m-1}}{}{}}
{\Res{V}{y_1 y_2 \cdots y_{m-1}}{}{}\cap (y_m)}
\longrightarrow
\frac{R/(y_m)}{\Tr{\HRes{I_t}{\y}{\p}{m-1}}{y_m}{j}{}}$$ is
injective. As we have $g_t \in \HRes{I_t}{\y}{\p}{m-1}$, it follows
that $$ \bar G_0 \in \frac{(\HRes{I_t}{\y}{\p}{m-1}+(y_m))_0}
{(y_m)}=\Tr{\HRes{I_t}{\y}{\p}{m-1}}{y_m}{0}{}, $$ i.e., $\varphi_0 (\bar
G_0)=0$, and therefore $\bar G_0=0$. Now we argue by iteration:
let $1\le j< p$, and assume we know $\bar G_0= \dots= \bar
G_{j-1}=0$. This means that $g_t \in (y_m, t^j)$, so
$G_j+\dots+G_{p-1}t^{p-1-j} \in
(\HRes{I_t}{\y}{\p}{m-1}+(y_m)):t^j$, which implies $\bar G_j \in
\Tr{\HRes{I_t}{\y}{\p}{m-1}}{y_m}{j+1}{},$ i.e., $\varphi_{j+1}
(\bar G_j)=0$, and therefore $\bar G_j=0$. 
\end{proof}

\begin{Cor}
\label{clm:formalah} Let $V \subset R$ be a $k$-linear subspace, and $I_t
\subset R_t$ an ideal such that $R_t/I_t$ is flat over 
  $k[[t]]$. Let $m\ge 0$ and let $\p=(p_1, \dots, p_m) \in \Z^m$ with
  $p_1 \ge p_2 \ge \dots \ge p_m \ge 1$, and $\y=(y_1, \dots, y_m) \in
  R^m$ be given sequences. 
  Suppose that the following two conditions are satisfied.
\begin{enumerate}
\item For $1 \le i \le m$, the canonical map
$$ \frac{\Res{V}{y_1 \cdots y_{i-1}}{}{}}{\Res{V}{y_1 \cdots
    y_{i-1}}{}{} \cap (y_i)} \longrightarrow
\frac{R/(y_i)}{\Tr{I_t}{\y}{\p}{i}}$$
is injective.
\item The canonical map
$$ \Res{V}{y_1 \cdots y_m}{}{} \longrightarrow
\frac{R}{\Res{I_t}{\y}{\p}{m}}$$
is injective.
\end{enumerate}
Then the canonical map
$\varphi_t: V \otimes k[[t]] \longrightarrow R_t/I_t$
is  injective.
\end{Cor}
\begin{proof}
  Let $f_t \in \Ker \varphi_t= V \otimes k[[t]] \cap I_t$. If
$f_0=0$, then we may replace $f_t$ by $f_t/t$ since $R_t/I_t$ is a flat
and hence torsion free $k[[t]]$-module. Thus we only have to prove
$f_0=0$. But proposition \ref{clm:formal} implies that 
$f_0= g y_1 \dots y_m$ for some $g \in\Res{V}{y_1 \cdots y_m}{}{}\cap
\Res{I_t}{\y}{\p}{m}$, so the second hypothesis 
gives $g_0=0$ and therefore $f_0=0$.
\end{proof}

\subsection{Preserving the number of conditions}
\label{sec:ncond}

We are mostly interested in flat families $R_t/I_t$ of finite length
(which define flat families of zeroschemes $Z_t$);
so it will be useful to consider the quantities
$\tr{I_t}{\y}{\p}{i}=\dim_k(R/\Tr{I_t}{\y}{\p}{i})$, 
$\hres{I_t}{\y}{\p}{i}=\dim_k(R_t/\HRes{I_t}{\y}{\p}{i})$, and
$\res{I_t}{\y}{\p}{i}=\dim_k(R/\Res{I_t}{\y}{\p}{i})$. In particular,
if $R_t/I_t$ is flat over $k[[t]]$ then
$\hres{I_t}{y}{\p}{0}=\infty$ and
$\res{I_t}{y}{\p}{0}=\dim_k R/I_0=\dim_{k((t))} (R_t/I_t)\otimes k((t)).$

Our aim is to obtain a linear system $\mathcal L$
which contains the limit of a family of linear systems $|L-Z_t|$
and, if possible, coincides with it. In the best cases, this will
serve to prove that general members of the family of
linear systems are regular, i.e., of dimension equal to $\dim |L| -
\len Z_t$, or $-1$ if this amount is negative. 

In the approach of section \ref{sec:formal},
$\mathcal L$ consists of the elements of $|L|$ that contain (a) the
divisors locally 
given by $y_1=0, \dots, y_m=0$ (containing $y_i=0$  accounts for
$\tr{I_t}{\y}{\p}{i}$ linear conditions) and (b) the residual
zeroscheme (which accounts for
$\res{I_t}{\y}{\p}{m}$ linear conditions). So if $\dim |L| -
\len Z_t \ge 0$, a requirement for the method to give
the desired result is that
\begin{equation}
  \label{eq:ncond}
  \res{I_t}{\y}{\p}{m}+\sum_{i=1}^m \tr{I_t}{\y}{\p}{i} 
=\res{I_t}{y}{\p}{0},
\end{equation}
and we now analyze when \eqref{eq:ncond} is satisfied. Note that
if it is not satisfied, the method can sometimes still be
applied to prove that a linear system of interest is empty (see
section \ref{sec:hhsquares}).

Given a sequence  $\p=(p_1, \dots,
p_m) \in \Z^m$ and an integer $q \le p_m$, define $\p-q=(p_1-q, \dots,
p_m-q)$.

\begin{Lem}
\label{clm:flatres}
  Let $I_t \subset R_t$ be an ideal such that $R_t/I_t$ is flat over
  $k[[t]]$. Let $m\ge 0$ and let $\p=(p_1, \dots, p_m) \in \Z^m$ with
  $p_1 \ge p_2 \ge \dots \ge p_m \ge 1$, and $\y=(y_1, \dots, y_m) \in
  R^m$ be given sequences. Then for every integer $q \le p_m$,
$$\HRes{I_t}{\y}{\p}{m}:t^q=\HRes{I_t}{\y}{\p-q}{m}.
$$
\end{Lem}
\begin{proof}
  Using \eqref{eq:altres}, it is easy to check that
  \begin{multline*}
   \HRes{I_t}{\y}{\p}{m}:t^q=
\left(I_t+(t^{p_1}, y_1t^{p_{2}},\ \dots\ ,
y_1 y_2 \cdots y_{m-1}t^{p_{m}})\right): (y_1 y_2 \cdots y_{m}t^q)= \\
=\left(I_t:t^q+(t^{p_1-q}, y_1t^{p_{2}-q},\ \dots\ ,
y_1 y_2 \cdots y_{m-1}t^{p_{m}-q})\right): (y_1 y_2 \cdots y_{m}).  
  \end{multline*}
The claim follows noting that $I_t:t^q=I_t$ (by flatness) 
and using \eqref{eq:altres} again. 
\end{proof}

It is well known that, for every ideal $I\subset R$, where $R$ is a
domain, and every $f\in R$, there is an exact sequence
\begin{equation}
\label{eq:ressq}
 0 \longrightarrow \frac{R}{I:f} \longrightarrow \frac{R}{I}
\longrightarrow \frac{R}{I+(f)} \longrightarrow 0,
\end{equation}
which we call
\emph{the residual exact sequence of $I$ with
respect to $f$}.
Given a sequence  $\p=(p_1, \dots,
p_m) \in \Z^m$ and two integers $j$ and $q$, with $1 \le j <m$ and
$q\le p_j$, let us denote $\p(q,j)=(p_1, p_2, \dots, p_{j-1},q)$.

\begin{Pro}
\label{clm:correct}
Let $I_t \subset R_t$ be an ideal such that $R_t/I_t$ is flat over
  $k[[t]]$ and $\dim_k R/I_0 < \infty$. Let $m\ge 0$ and let $\p=(p_1, \dots, p_m) \in \Z^m$ with
  $p_1 \ge p_2 \ge \dots \ge p_m \ge 1$, and $\y=(y_1, \dots, y_m) \in
  R^m$ be given sequences. Then
\begin{enumerate}
  \item for every $i=1, \dots, m$, \begin{multline*}
  \res{I_t}{\y}{\p}{i}+\sum_{j=1}^{i}\tr{I_t}{\y}{\p}{j}=
\res{I_t}{\y}{\p}{0}-\\
-\sum_{j=1}^{i} 
\dim \frac{\HRes{I_t}{\y}{\p-1}{j-1}+(t^{p_j-1},y_j)}{\HRes{I_t}{\y}{\p}{j-1}+(t^{p_j-1},y_j)},
\end{multline*}
  \item for every $j=1, \dots, i$,
$$
\dim
\frac{\HRes{I_t}{\y}{\p-1}{j-1}+(t^{p_j-1},y_j)}
{\HRes{I_t}{\y}{\p}{j-1}+(t^{p_j-1},y_j)} 
=\sum_{q=1}^{p_j-1} 
\left( \tr{I_t}{\y}{\p(q,j)}{j}-
  \tr{I_t}{\y}{(\p-1)(q,j)}{j} \right). 
$$
\end{enumerate}
\end{Pro}

The first claim of proposition \ref{clm:correct} gives the amount by which the
higher traces and residuals of $I$ with respect to $\y$ fail to preserve
the number of conditions imposed to the linear system.
The second shows that this amount can be exactly computed
whenever we can compute the colengths of the traces (even
if the residuals are unknown).

\begin{proof}
Applying the residual exact sequence of 
$\HRes{I_t}{\y}{\p}{i-1}+(t^{p_i})\subset R_t$ with respect to $y_i\in
R_t$ gives
$$\hres{I_t}{\y}{\p}{i}=\dim \frac{R_t}{\HRes{I_t}{\y}{\p}{i-1}+(t^{p_i})} - 
\dim \frac{R_t}{\HRes{I_t}{\y}{\p}{i-1}+(t^{p_i},y_i)}.$$

The two terms on the right can be evaluated by means of residual
exact sequences. Indeed, the residual exact sequence of 
$\HRes{I_t}{\y}{\p}{i-1}+(t^{p_i})$ with respect to $t^{p_i-1}$, 
together with lemma  \ref{clm:flatres},
gives 
$$\dim \frac{R_t}{\HRes{I_t}{\y}{\p}{i-1}+(t^{p_i})}=
\dim \frac{R_t}{\HRes{I_t}{\y}{\p}{i-1}+(t^{p_i-1})}+
\res{I_t}{\y}{\p-p_i+1}{i-1},$$
and recursively applying this last equality,
$$\dim \frac{R_t}{\HRes{I_t}{\y}{\p}{i-1}+(t^{p_i})}=
\sum_{q=1}^{p_i}\res{I_t}{\y}{\p-q+1}{i-1}.$$

 On the other hand, applying the residual exact sequence of 
 $\HRes{I_t}{\y}{\p}{i-1}+(t^{p_i},y_i)$ with respect to $t^{p_i-1}$
 we get
 \begin{equation}
   \label{eq:3}
   \dim \frac{R_t}{\HRes{I_t}{\y}{\p}{i-1}+(t^{p_i},y_i)}=
 \dim \frac{R_t}{\HRes{I_t}{\y}{\p}{i-1}+(t^{p_i-1},y_i)}+
 \tr{I_t}{\y}{\p}{i},
 \end{equation}
and so
$$\hres{I_t}{\y}{\p}{i}=\sum_{q=1}^{p_i}\res{I_t}{\y}{\p-q+1}{i-1}-
\tr{I_t}{\y}{\p}{i}- \dim \frac{R_t}{\HRes{I_t}{\y}{\p}{i-1}+(t^{p_i-1},y_i)}.$$


Now apply the residual exact sequence of $\HRes{I_t}{\y}{\p}{i}$ with
respect to $t$ and lemma \ref{clm:flatres}
to obtain that
$\res{I_t}{\y}{\p}{i}=\hres{I_t}{\y}{\p}{i}-\hres{I_t}{\y}{\p-1}{i}.$

Putting together everything we have so far, it follows that
\begin{multline*}
  \res{I_t}{\y}{\p}{i}=\res{I_t}{\y}{\p}{i-1}-\tr{I_t}{\y}{\p}{i}-\\
-\left(\dim \frac{R_t}{\HRes{I_t}{\y}{\p}{i-1}+(t^{p_i-1},y_i)}-
\dim \frac{R_t}{\HRes{I_t}{\y}{\p-1}{i-1}+(t^{p_i-1},y_i)}\right),
\end{multline*}
which recursively applied yields
the first claim. The second follows by
applying (\ref{eq:3}) recursively.
\end{proof}

One implication of proposition \ref{clm:correct} is that the number of
conditions is preserved whenever the integers in the sequence $\p$
decrease ``fast enough'', which will give the second part of theorem
\ref{clm:global}. We prove this next: 

\begin{Cor}
\label{clm:decreasp}
  Assume that 
$$\res{I_t}{\y}{\p}{i}+\sum_{j=1}^{i}\tr{I_t}{\y}{\p}{j}=
\res{I_t}{\y}{\p}{0}.$$
Then there exists an integer $q_i$, $p_i \ge q_i \ge 0$ such that
$$\res{I_t}{\y}{\p}{i+1}+\sum_{j=1}^{i+1}\tr{I_t}{\y}{\p}{j}=
\res{I_t}{\y}{\p}{0} \text{ if and only if } p_{i+1}\le q_i.$$
\end{Cor}

\begin{proof}
  The hypothesis, together with proposition \ref{clm:correct}, tell us that
\begin{multline*}
  \res{I_t}{\y}{\p}{i+1}+\sum_{j=1}^{i+1}\tr{I_t}{\y}{\p}{j}=\\
=\res{I_t}{\y}{\p}{0}
-\dim \frac{\HRes{I_t}{\y}{\p-1}{i}+(t^{p_{i+1}-1},y_{i+1})}
{\HRes{I_t}{\y}{\p}{i}+(t^{p_{i+1}-1},y_{i+1})}=\\
=\res{I_t}{\y}{\p}{0}-\sum_{q=1}^{p_i-1} 
\left( \tr{I_t}{\y}{\p(q)}{i}- \tr{I_t}{\y}{(\p-1)(q)}{i} \right) .
\end{multline*}
Now, for all $q$ it is easy to see that 
$\Tr{I_t}{\y}{\p(q)}{i} \subseteq \Tr{I_t}{\y}{(\p-1)(q)}{i}$, 
so clearly $q_i=\min \left\{ q \left|
  \Tr{I_t}{\y}{\p(q)}{i} \neq \Tr{I_t}{\y}{(\p-1)(q)}{i} \right.\right\}-1$
satisfies the claim.
\end{proof}

\begin{Exa}
\label{clm:1layer}
   It follows from corollary \ref{clm:decreasp} that non iterated traces
   and residuals (i.e., when $m=1$, which is the case used in \cite{AH00}) 
   always preserve the number of conditions. On the other hand, it
   follows from \cite{Eva??}, \cite{Eva07} that for monomial
   schemes approaching a unique divisor $y=0$ vertically,
   $q_i=p_i-1$. The simplest examples in which $p_{i+1}< p_i$ for all
   $i$ but the number of conditions is not preserved involve
   monomial schemes moving non-vertically. Let $f=x+y+t\in
   R_t=k[[x,y]][[t]]$, $I=(f,x^2)^4$,
   $\y=(y,y,x)$, and $\p=(8,7,p_3)$. It is not hard to compute $q_2=5$;
   therefore the number of conditions is not preserved if $\p=(8,7,6)$.
\end{Exa}

\begin{proof}[Proof of theorem \ref{clm:global}]
Let $R$ be the local ring of $X$ at the support point of $Y_{t_0}$ (or its
completion), $y_i$ a local
equation of $D$, $V$ the image of the natural map $\O_X(L)\rightarrow
R$ and $k[[t]]$ the completion of the local ring of $C$ at $t_0$. Then
proposition \ref{clm:formal} gives the first part of the statement.

The second part of the statement follows from \ref{clm:decreasp}. Indeed,
assume $\rho_i^D$ is bijective for all $i$. If 
$H^0(X,\O_X(L-(i-1)D))\rightarrow H^0(D,\O_D(L-(i-1)D))$ has maximal 
rank for all $i$ and is injective for some $i$, then 
$\edim |L-Z_t|= \edim |(L-mD)-\Res{I_t}{\y}{\p}{m}|=\dim |L-mD|=0,$
and if $H^0(X,\O_X(L-(i-1)D))\rightarrow H^0(D,\O_D(L-(i-1)D))$
is surjective for all $i$, then 
$\dim |L-iD|=\dim |L-(i-1)D|-\len \TrZ{Z_t}{D}{\p}{i}$. 
If $p_{i+1}\le q_i$ with $q_i$ as in \ref{clm:decreasp}, then 
$$\ResZ{Z_t}{D}{\p}{i}+\sum_{i=1}^{m}\tr{Z_t}{D}{\p}{i}=
\res{Z_t}{D}{\p}{0}$$
so the claim follows.
\end{proof}

\subsection{A computation for squares}
\label{sec:sqcomp}

Note that even if the number of conditions is not preserved, the
method can still be useful to show that a linear system is empty.
We illustrate this with a computation in $R_t=k[[x,y]][[t]]$
(proposition \ref{clm:sqodd})
needed for the proof of theorem \ref{clm:squares}.
To simplify, consider the
leading, or dominant, terms of series  in the local ring $k[[x,y]][[t]]$,
with respect to some regular system of parameters of the form
$\{x,f,t\}$. Given  
$g=\sum a_{ijk}x^if^jt^k \in k[[x,f]][[t]]=k[[x,y]][[t]]$, we set 
$\ord_t(g)=\min \{k| \exists i,j; \, a_{ijk}\ne 0 \}$ and define the
dominant part of $g$ as
$\lead{g}=\sum a_{i,j,\ord_t(g)}x^if^jt^{\ord_t(g)}$. An
ideal $I_t \subset k[[x,y,t]]$ determines  its  ideal of
dominant terms $\lead{I}_t=\left(\lead{g}\right)_{g \in I_t}\subset
k[[x,y,t]]$. Observe that $(\lead{I}_t)_0=I_0$.

\begin{Lem}
  Let $f=y-t\in k[[x,y]][[t]]$, $I_t\subset k[[x,y,t]]$ and $p
  \in \Z$, $p \ge 1$. Assume that $\lead{I}_t \subset (x,f)^m+(t^p)$. Then
  $\lead{\left((I_t+(t^{p})):y\right)} \subset (x,f)^{m-1}+(t^p)$
\end{Lem}

\begin{proof}
  Let $g =\sum a_{ijk}x^if^jt^k \in k[[x,f]][[t]]$ be such that
$yg \in I_t+(t^{p})$. We want to see that
$\lead{g}\in(x,f)^{m-1}+(t^p)$.
But $yg=(f+t)g=\sum
(a_{ijk}+a_{i,j,k-1})x^if^jt^k$ (where we set $a_{i,j,-1}=0$ for all
$i, j$) so $\lead{(yg)}=f(\lead{g})$. On the other hand, $\lead{(yg)}
\in \lead{(I_t+(t^p))}=\lead{I}_t+(t^p)\subset (x,f)^m+(t^p)$. Both
things together tell us that $\lead{g}\in \left((x,f)^{m}+(t^p)
\right):f =(x,f)^{m-1}+(t^p)$, as claimed.
\end{proof}

\begin{Pro}
\label{clm:traces}
  Let $e$, $p$, be positive integers with $e +1\ge p$, and let $I_t
  =(x,y-t) \subset k[[x,y,t]]$. Then
  \begin{enumerate}
  \item For every positive integer $i$, 
$\tr{I^e_t}{y}{p}{i} \ge e+2-p-i,$ and
  \item Assume $k$ has characteristic zero. Then $\Res{I^e_t}{y}{p}{e} \subset \M^{\lfloor{p \over2}\rfloor}$.
  \end{enumerate}
\end{Pro}

\begin{proof}
  Using the previous lemma $i-1$ times, it follows that
  $\HRes{I^e_t}{y}{p+1}{i-1}+(t)\subset I^{e-i+1}_t+(t)$. On the other
  hand, it is easy to see (and is proved as part of proposition 8.1 in
  \cite{AH00}) that $\tr{I^{e-i+1}_t}{y}{p}{} = e+2-p-i,$
  whence the first claim.

  Because of (\ref{eq:altres}), what remains to prove is
$$ 
\left(\left(I_t^e+(t^{p})\right): y^{e}\right)_0 \subset
(x,y)^{\lfloor{p \over2}\rfloor}.$$ 
Define again $f=y-t$ and consider the
automorphism $\varphi$ of  $R_t$ defined by 
  $\varphi(x)=x$, $\varphi(y)=f$, $\varphi(t)=t$. It is a
  $k[[t]]$-automorphism (it leaves $k[[t]]$ fixed) and 
for every $g=\sum a_{ijk}x^iy^jt^k \in
k[[x,y,t]]$, $\varphi (g)_0=\sum a_{ij0}x^iy^j=g_0$. So the claim is
equivalent to
$$ 
\left(\left(\varphi^{-1}(I_t^e)+(t^{p})\right):
  \varphi^{-1}(y)^{e}\right)_0=
\left(\left((x,y)^e+(t^{p})\right): (y+t)^{e}\right)_0 \subset
(x,y)^{\lfloor{p \over2}\rfloor},
$$
i.e., if 
$g(y+t)^e\in (x,y)^e+(t^{p})$ then we need to prove
$g_0 \in (x,y))^{\lfloor{p \over2}\rfloor}$. But now 
the ideal $\varphi^{-1}(I_t^e)+(t^p)=(x,f)^e+(t^{p})$ is monomial,
and so $h=\sum b_{ijk}x^iy^jt^k \in (x,y)^e+(t^{p})$ if and only if
$b_{ijk}=0$ for all $k <p$ and $i+j < e$.

Let now $g=\sum a_{ijk}x^iy^jt^k$ and assume that 
$h=g(y+t)^e=\sum b_{ijk}x^iy^jt^k \in (x,y)^e+(t^{p})$. By definition,
$$h=\sum a_{ijk} \sum_{\ell=0}^e {e \choose \ell}
x^{i}e^{\ell+j}t^{k+m-\ell},$$
so 
$$b_{ijk}=\sum_{\ell=\ell_0}^{\ell_1}  {e \choose \ell} 
a_{i,j-\ell,k+\ell-e},$$ 
where $\ell_0=\max \{e-k,0\}$ and  $\ell_1=\max \{j,e\}$ .
The condition $h \in I_t^e$ thus translates into the linear equations
$$ \sum_{\ell=e-k}^{j}  {e \choose \ell} a_{i,j-\ell,k+\ell-e}, \qquad
0\le i,j,k; \ k <p; \ i+j \le e-1.$$ 
Some among these equations involve the same set of coefficients;
namely, for each fixed $i$ and $r=j+k-e$ satisfying $0\le i \le e-1$ and $0\le
r < p-i-1$ we have obtained a system of linear equations
\begin{equation}
  \label{eq:2}
  \begin{pmatrix}
    {e \choose {e-p+1}} & {e \choose {e-p+2}} & \cdots & {e \choose r}\\ 
    {e \choose {e-p+2}} & {e \choose {e-p+3}} & \cdots & {e \choose {r+1}}\\
       \vdots     &     \vdots    & \ddots &   \vdots  \\
 {e \choose {e-1-i-r}}&{e\choose{e-i-r}}&\cdots&{e\choose{e-1-i}}
  \end{pmatrix} 
  \begin{pmatrix}
    a_{i,r,0} \\ a_{i,r-1,1} \\ \vdots \\ a_{i,0,r}
  \end{pmatrix}
=0,
\end{equation}
which, if $e-i-1-r-(e-p+1)\ge r$, admits only the trivial solution
because the matrix on the left has nonzero determinant
(see lemma \ref{lem:det} below). In particular, if $h \in I_t^e$ then for
every $(i,r)$ with  $i+r< \lfloor{p \over2}\rfloor$ 
(which trivially implies $e-i-1-r-(e-p+1)\ge r$) we obtain $a_{i,r,0}=0$,
which means $g_0 \in (x,y)^{\lfloor{e \over2}\rfloor}$, and the second
claim follows.
\end{proof}

\begin{Lem}
\label{lem:det}
  For every triple of integers $e\ge r \ge n\ge 1$, the
  following symmetric $(n+1)\times(n+1)$ matrix is invertible.
\begin{equation*}
  H_{r,n}(e)=\begin{pmatrix}
    {e \choose r-n} & {e \choose r-n+1} & \cdots & {e \choose r}\\ 
    {e \choose r-n+1} & {e \choose r-n+2} & \cdots & {e \choose {r+1}}\\
       \vdots     &     \vdots    & \ddots &   \vdots  \\
 {e \choose {r}}&{e\choose{r+1}}&\cdots&{e\choose{r+n}}
  \end{pmatrix} 
\end{equation*}
\end{Lem}

Similar matrices are known in the litterature, e.g., in
\cite[p. 94]{ACGH85}, \cite{Rad00}, \cite{Jun03}.

\begin{proof}
  For $r,n$ fixed, $\det H_{r,n}(e)\in \Q[e]$ is a
  polynomial of degree (at most) $r(n+1)$, since the entry in the
  $(i,j)$ position, $0\le i,j\le n$, is
$$
\binom{e}{r-n+i+j}=\frac{\prod_{a=0}^{r-n+i+j-1}(e-a)}{(r-n+i+j)!},
$$
a polynomial of degree $r-n+i+j$. Moreover, every element of the
$i$th row is divisible by $\binom{e}{r-n+i}$, and therefore $\det
H_{r,n}(e)$ is divisible by
$P(e)=\prod_{i=0}^n \binom{e}{r-n+i}.$
On the other hand, by using the identities
$$
\binom{e}{a}+\binom{e}{a+1}=\binom{e+1}{a+1},
$$
a few elementary operations on rows show that
\begin{equation*}
  \det H_{r,n}(e)=\det \begin{pmatrix}
    {e \choose r-n} & {e \choose r-n+1} & \cdots & {e \choose r}\\ 
    {e+1 \choose r-n+1} & {e+1 \choose r-n+2} & \cdots & {e+1 \choose {r+1}}\\
       \vdots     &     \vdots    & \ddots &   \vdots  \\
 {e+n \choose {r}}&{e+n\choose{r+1}}&\cdots&{e+n\choose{r+n}}
  \end{pmatrix}.
\end{equation*}
Taking again common divisors to elements in each row one gets that 
$\det H_{r,n}(e)$ is divisible by
$Q(e)=\prod_{i=0}^n \binom{e+i}{r-n+i}.$
Therefore $\det H_{r,n}(e)$ is divisible by the $\lcm$ of $P$ and $Q$,
 which has degree
$r(n+1)$ and all its roots in the set $\{-n, -n+1, \dots, r-1\}$. It
follows that all roots of $\det H_{r,n}(e)$ are strictly less than $r$, and
thus for $e\ge r$, this determinant does not vanish.
\end{proof}

\begin{Pro}
\label{clm:sqodd}
  Assume that the characteristic of the base field $k$ is zero.
  Let $n=s^2$ be an odd square, and let $e>s/2$ be an integer. Let
  $a=se+(s-5)/2$. Given a set of $n$ distinct points in general position in
  $\P^2$, there are no curves of degree $a$ with multiplicity at least
  $e$ at each point.
\end{Pro}

\begin{proof}
Consider an irreducible smooth curve $C$ of
degree $s$ (and genus $g=(s-1)(s-2)/2$), and let $p_1, \dots,
p_{s^2-1}$ be general points of $C$, 
whereas $p_{s^2}$ is a general point of $\P^2$. Denote by $Z$ the
union of these points taken with multiplicity $e$. The restriction of
$\I_{Z}(a)$ to $C$ is an invertible sheaf of degree
$d=as-(s^2-1)m=m-s(s-5)/2$ which, by the genericity of the choice of
the $s^2-1>g$ points, is general among those of its degree. 

If $d<g$ then this invertible sheaf has no nonzero global sections,
i.e., the curve $C$ is a fixed part of the linear system
$H^{0}(\I_{Z}(a))$. The residual linear system is formed by curves of 
degree $a-s= s(e-1)+(s-5)/2$ which contain the scheme $Z'$ consisting
of the points have $p_1, \dots, p_{s^2-1}$ with multiplicity $e-1$ 
and $p_{s^2}$ with multiplicity $e$. But then the restriction of
$\I_{Z'}(a-s)$ to $C$ is an invertible sheaf of degree $d'=d-1<g$, and is
still general among those of its degree, so $C$ is again a fixed part
of the linear system. Iterating this process, we see that $C$ is
contained exactly $e$ times in the curves of the linear system
$H^{0}(\I_{Z}(a))$, and the residual linear system consists of curves
of degree $a-sm=(s-5)/2$ with a point of multiplicity $e>s/2$ so it is
empty as claimed. 

So assume $d\ge g$ and let $p=e+g-d=s+2$ (so we trivially have $e+1\ge
p \ge 1$). Now let $p_{s^2}$ tend
to $C$ transversely. In other words, choose a (general) point $q\in
C$, let $x,y$ be local parameters at $q$ such that $y=0$ is a
local equation for $C$, and let $q_t=(0,-t)$. We want to see that the
limit of the linear systems formed by curves of degree $a$ with
multiplicity $e$ at $p_1, \dots, p_{s^2-1}$ and at $q_t$ when
$t\mapsto 0$ is empty. Applying the first claim of \ref{clm:traces} and
\ref{clm:formal}, it follows that the limit linear system consists of
$C$ counted $e$ times plus a moving part, consisting of curves of
degree $a-es= (s-5)/2$ going through the zeroscheme defined by 
$\Res{I^e_t}{y}{p+1}{e}$. But by the second claim of \ref{clm:traces}
this zeroscheme contains the point $q$ counted $p/2>s/2$ times,
and we are done.
\end{proof}

\subsection{Homogeneous ideals in power series rings}
\label{sec:homog}


Throughout this section we assume that 
$R= k[[x_1, \dots, x_r]]$ is a power series ring, so
both $R$ and $R_t$ are regular local rings, whose maximal ideals are 
$\mathfrak{m}=(x_1, \dots, x_r)$ and $\mathfrak{m}_t=(x_1, \dots,
x_r,t)$ respectively, and come endowed with a natural
$\mathfrak{m}$-adic valuation $v$. 

Let $I_t \subset R_t$ be a
homogeneous ideal, and $\y=(y_1, \dots, y_m) \in R_t^m$ a sequence of
homogeneous polynomials. All higher traces and
residuals defined above are then homogeneous. Moreover, it is
easy to see that if $I$ and $J$ are homogeneous ideals and $f$ is
a homogeneous polynomial, then for every integer $e$
\begin{equation}
  \label{eq:reshom}
I+\mathfrak{m}^e= J+\mathfrak{m}^e \ \Rightarrow
I:f+\mathfrak{m}^{e-v(f)}= J:f+\mathfrak{m}^{e-v(f)}.
\end{equation}
Thus higher traces and residuals do not
differ from ordinary traces and residuals (up to a
finite order that can be computed). So we will substitute
one for the other to avoid too cumbersome computations (with due cautions,
essentially contained in the following proposition).

\begin{Pro} \label{clm:cutideal1}
  Let $I_t \subset R_t= k[[x_1, \dots, x_r,t]]$ be a homogeneous ideal.
  Let $m\ge 0$ and let
  $\p=(p_1, \dots, p_m) \in \Z^m$ and $\y=(y_1, \dots, y_m) \in
  R^m$ be given sequences, with $y_i$
  homogeneous for all $i$, such that  $p_i - p_{i+1}\ge v(y_i)$ for
  $i=1, \dots, m-1$ and $p_m\ge v(y_m)$. For each $j=1, \dots, m$,
  define $\y_j=(y_j, \dots,  y_m)$ and $\p_j=(p_j, \dots,
  p_m)$; for $1\le j\le i \le m$, let
  $V_j^i=\sum_{\ell=j}^{i}v(y_\ell)$. Then, 
  \begin{enumerate}
  \item \label{ordhres} for every $1\le j\le i \le m$,
$$\HRes{I_t}{\y}{\p}{i}+\mathfrak{m}_t^{p_j - V_j^i}
      =\HRes{I_t:(y_1\dots y_{j})}{\y_j}{\p_j}{i-j}+
      \mathfrak{m}_t^{p_j -V_j^i},$$ 
  \item \label{ordtr} for every $1\le j<i \le m,$
 $$\Tr{I_t}{\y}{\p}{i}+\mathfrak{m}^{p_{j} -
      p_{i}+1 - V_{j}^{i-1}}=\Tr{I_t:(y_1\dots
      y_{j})}{\y_{j}}{\p_j}{i-j}+ \mathfrak{m}^{p_{j} -
      p_{i}+1 -V_{j}^{i-1}}.$$
  \end{enumerate}
\end{Pro}

\begin{proof}
Observe that the first claim, together with
\eqref{eq:reshom}, implies the second. 
Again, using \eqref{eq:reshom} $i-j$ times we see that the first claim
will follow from
$$\HRes{I_t}{\y}{\p}{j}+\mathfrak{m}_t^{p_j -
      v(y_j)}=I_t:(y_1\dots y_{j})+
      \mathfrak{m}_t^{p_j - v(y_j)},$$
which we prove by induction on $j$.
For $j=0$, there is nothing to prove. For $j>0$, since
 $p_j \le p_{j-1}-v(y_{j-1})$ the induction hypothesis gives
$$\HRes{I_t}{\y}{\p}{j-1}+(t^{p_j})+\mathfrak{m}_t^{p_{j}}=
\HRes{I_t}{\y}{\p}{j-1}+\mathfrak{m}_t^{p_{j}}=I_t:(y_1\dots y_{j-1})+
      \mathfrak{m}_t^{p_j},$$
which using (\ref{eq:reshom}) and the definition of
$\HRes{I_t}{\y}{\p}{j}$ finishes the proof.
\end{proof}

The following technical lemma, which is a rather straightforward
application of proposition \ref{clm:cutideal1}, is the key to show
that the number of conditions is preserved in all the applications
contained in this paper. Given an ideal $I\subset R_t= k[[x_1,
\dots, x_r,t]]$ and sequences  $\p=(p_1, \dots, p_m) \in \Z^m$ and
$\y=(y_1, \dots, y_m) \in R^m$, denote
$$
\val{I}{\y}{\p}{i} = \max_{1\le q\le p_i} \{\tr{I}{\y}{\p(q,i)}{i}+q-p_i\}.
$$

\begin{Lem} \label{clm:cutideal2} 
  Let $I_t \subset R_t= k[[x_1, \dots, x_r,t]]$ be a homogeneous ideal
  such that $R_t/I_t$ is flat over $k[[t]]$. Let $m\ge 0$ and let
  $\p=(p_1, \dots, p_m) \in \Z^m$ and $\y=(y_1, \dots, y_m) \in
  R^m$ be given sequences, with $y_i$
  homogeneous for all $i$, such that  $p_i - p_{i+1}\ge v(y_i)$ for
  $i=1, \dots, m-1$ and $p_m\ge v(y_m)$. For each $j=1, \dots, m$,
  define $\y_j=(y_j, \dots,  y_m)$ and $\p_j=(p_j, \dots,
  p_m)$; for $1\le j\le i \le m$, let
  $V_j^i=\sum_{\ell=j}^{i}v(y_\ell)$.
  Assume that for every $i>1$ there is $j<i$ with
  \begin{enumerate}
  \item \label{condp} 
        $p_{j} - p_{i} \ge V_{j}^{i-1}+\val{I_t:(y_1
        \dots y_{j})}{\y_j}{\p_i(p_i-1)}{i-j}-2,$
    \item $\res{I_t:(y_1\dots
        y_{j})}{\y_j}{\p_j}{i-j}+\sum_{\ell=i+1}^j \tr{I_t:(y_1
        \dots y_{j})}{\y_j}{\p_j}{\ell}  
=\newline =\res{I_t:(y_1\dots y_{j})}{\y_j}{\p_j}{0}.$
  \end{enumerate}
Then $$  \res{I_t}{\y}{\p}{m}+\sum_{i=1}^m \tr{I_t}{\y}{\p}{i} 
=\res{I_t}{\y}{\p}{0}.$$
\end{Lem}

\begin{proof}
Due to proposition \ref{clm:correct}, what is needed to show is that for
each $i=1,\dots,m$ and $q=1,\dots, p_i-1$, 
$\tr{I_t}{\y}{\p(q)}{i} = \tr{I_t}{\y}{(\p-1)(q)}{i}.$
By hypothesis there is $j<i$ with $p_{j} -
q-\sum_{\ell=j}^{i-1}v(y_\ell) +1 \ge 
\tr{I_t:(y_1\dots y_{j})}{\y_j}{\p_j(q)}{i-j}$ for all $1\le q \le
p_i-1$, and therefore $\mathfrak{m}^{p_{i-1} - q
  -\sum_{\ell=j}^{i-1}v(y_{\ell})+1}$ is contained in $\Tr{I_t:(y_1\dots
  y_{j})}{\y_j}{\p_j(q)}{i-j}$. We also have
\begin{align*}
  \Tr{I_t:(y_1\dots
  y_{j})}{\y_j}{\p_j(q)}{i-j} & \subset \Tr{I_t:(y_1\dots
  y_{j})}{\y_j}{(\p_j-1)(q)}{i-j} \\
\Tr{I_t:(y_1 \dots y_{j})}{\y_i}{\p_j(q)}{i-j} & \subset
  \Tr{I_t}{\y}{\p(q)}{i} \\
\Tr{I_t:(y_1 \dots y_{j})}{\y_i}{(\p_j-1)(q)}{i-j} & \subset
  \Tr{I_t}{\y}{(\p-1)(q)}{i},
\end{align*}
so in particular $\mathfrak{m}^{p_{i-1} -
  q-V_j^{i-1}+1}$ is contained in the four ideals
  involved.  Using the second part of proposition \ref{clm:cutideal1} we get
\begin{multline*}
\Tr{I_t}{\y}{\p(q)}{i}=\Tr{I_t}{\y}{\p(q)}{i}+\mathfrak{m}^{p_{i-1} -
  q-V_j^{i-1}+1} = \\
=\Tr{I_t:(y_1\dots y_{j})}{\y_j}{\p_j(q)}{i-j}+\mathfrak{m}^{p_{i-1} -
  q-V_j^{i-1}+1} = \\
=\Tr{I_t:(y_1\dots y_{j})}{\y_j}{(\p_j-1)(q)}{i-j}+\mathfrak{m}^{p_{i-1} -
  q-V_j^{i-1}+1} = \\
=\Tr{I_t}{\y}{(\p-1)(q)}{i}+\mathfrak{m}^{p_{i-1} -
  q-V_j^{i-1}+1} =\Tr{I_t}{\y}{(\p-1)(q)}{i},
\end{multline*}
the third equality being a consequence of the second hypothesis and
proposition \ref{clm:correct}. Now the result follows by \ref{clm:correct}.
\end{proof}

\subsection{Staircases and monomial ideals}
\label{sec:monomial}

We now specialize to the two-dimensional case, so let $R=k[[x,y]]$.
Given a \emph{staircase} $E \subset \Z_{\ge 0}^2$, i.e., a subset
satisfying $E + \Z_{\ge 0}^2\subset E$, and a system of
parameters $f,g \in R \cong k[[x,y]]$, we denote
$$I_{E,f,g}=\left(f^{e_1}g^{e_2}\right)_
{(e_1,e_2)\in E}.$$
  If $E\subset \Z_{\ge 0}^2$ is a staircase,  the
  \emph{length} of its $i$th stair is $\ell_E(i)=\min \{e\ |\ (e,i)\in
  E\}$, and the \emph{height} of its $i$th slice is $h_E(i)=\min
  \{e\ |\ (i,e)\in E\}$. We use the first difference of $\ell$
  as well: $\hat \ell_E(i)=\ell_E(i)-\ell_E(i+1)$. When there are no
  steps of height $>1$, i.e., if $h_E(i)\le
  h_E(i+1)+1$ for all $i$, we say that $E$ is \emph{gentle}. We
  also define the total length and height of $E$ as 
  $\ell(E)=\ell_E(0)$ and $h(E)=h_E(0)$, and the minimal length
  $\lmin(E)=\min\{\hat\ell_E(i)\ |\ 0\le i<h(E)-1\}$ (for technical reasons that
  will become apparent in forthcoming sections, the latter does not take into
  account the length of the top stair).

\begin{Lem}
\label{clm:indepx}
  For every staircase $E$ with finite complement, and every system of
  parameters $f,g \in R \cong k[[x,y]]$,
\begin{enumerate}
\item $I_{E,f,g}$ is $\M$-primary, and has colength $\#(\Z_{\ge 0}^2
\setminus E)$, and
\item $I_{E,f,g}$ depends only on finite jets of $f$ and $g$, i.e.,
  there exist integers $a=a(E)$ and $b=b(E)$ such that $f_1-f_2\in
  \M^a$, $g_1-g_2   \in \M^b$ imply
  $I_{E,f_{1},g_{1}}=I_{E,f_{2},g_{2}}$.
\item if $E$ is gentle then $I_{E,f,g}$ does not depend on $f$, i.e.,
  $I_{E,f_{1},g}=I_{E,f_{2},g}$ whenever
  $(f_{1},g)=(f_{2},g)=\M$. In such a case we denote
  $I_{E,g}=I_{E,f_{1},g}$. 
\end{enumerate}
\end{Lem}

\begin{proof}
Because $E$ has finite complement, it follows that for suitable $e_1,
  e_2$, $f^{e_1} \in I_1$ and  $g^{e_2} \in I_{E,f,g}$, so 
$$\M^{e_1+e_2}=(f,g)^{e_1+e_2}\subset I_{E,f,g}$$ and $I_{E,f,g}$
is $\M$-primary. Observe that $e_1, e_2$ depend only on $E$, not on
$f$ or $g$. The colength follows from the well known fact that
the classes modulo $I_{E,f,g}$ of the monomials $f^{e_1}g^{e_2}$ with
$(e_1,e_2)$ \emph{not} in $E$ form a basis of $R/I_{E,f,g}$.

  For the second claim, by symmetry, it is enough to prove 
  that $I_{E,f,g}$ depends only on a finite jet of $f$. We have just
  seen that there is a fixed integer $a$ such that $\M^a \subset I_{E,f,g}$
  for every choice of $f$, and we want to prove that
  given $f_1, f_2 \in R$ with $f_1-f_2 \in \M^a$,
  $I_{E,f_{1},g}=I_{E,f_{2},g}$. Again by symmetry it 
  will be enough to show that for every $(e_1,e_2)\in E$, 
  $f^{e_1}_{1}g^{e_2} \in I_{E,f_{2},g}$. This follows from
$$f^{e_1}_{1}g^{e_2}-f^{e_1}_{2}g^{e_2}=g^{e_2}(f_{1}^{e_1}-f_{2}^{e_1})\in
g^{e_2}(f_{1}-f_{2}). 
$$

Finally for the third claim, and by symmetry again, we have to see
that given $f_1, f_2$ such that $(f_1,g)=(f_2,g)=\M$, for every
$(e_1,e_2)\in E$ one has $f^{e_1}_1g^{e_2} \in I_{E,f_2,g}$. But,
because the staircase is gentle, 
  it follows that if $(e_1,e_2)\in E$, then for every integer $0\le k
  \le e_1$, $(e_1-k, e_2+k)\in E$, and therefore
  $g^{e_2}(f_2,g)^{e_1}\subset I_{E,f_2,g}$. Now since
  both $(f_1,g)$ and $(f_2,g)$ are systems of parameters, it follows
  that $f_1 \in(f_2,g)$, and 
  $f^{e_1}_1g^{e_2}\in g^{e_2}(f_2,g)^{e_1}\subset I_{E,f_2,g}.$
\end{proof} 

We are interested in a specific kind of families of translated
monomial ideals. Fix $f=x+y+t \in R_t=k[[x,y,t]]$. For every staircase
$E$, define
$$
I_E=I_{E,x,f}=\left(x^{e_1}f^{e_2}\right)_
{(e_1,e_2)\in E}.
$$
The sequences $\y$ of interest will have a
fixed form as well, namely
$$
\y_{\m}=(y,x,\dots,x,y,x,\dots,x,y,x,\dots)
$$
where $\m=(m_1,\dots,m_\mu) \in \Z_{> 0}^\mu$, 
$y_i=y$ for $i =1, 1+m_1, \dots, 1+m_1+ \dots +m_{\mu-1}$ and 
$y_i=x$ for all other $i \le \sum m_i$. In other words, $\y_\m$ is
the concatenation of $\mu$ sequences of the form
$(y,x,\dots,x)$ of lengths $m_1,\dots,m_\mu$.

The properties of staircase ideals 
with respect to higher order traces and residuals have been
extensively studied by \'Evain \cite{Eva??}, \cite{Eva07}, in the
particular case of vertical translations (roughly
speaking, using $f=y+t$). His results show that
traces can be computed from slices of the staircase, and residuals are
obtained by deleting the same slices. This is not always the case for
non-vertical translations like the ones just defined, as showed by
example \ref{clm:1layer}; the key lemma \ref{clm:cutideal2} will
show that under suitable numerical conditions \'Evain's computations
do hold in our setting as well.

For convenience, we introduce a function $\sigma_\m$ to count the number
of $x$ appearing in $\y_\m$ up to the $i$th position, and horizontal
translation of staircases.
\begin{align*}
  \sigma_\m(i)&=i-1-\max\left\{k \left| \sum_{j=1}^km_j\le
  i-1\right\}\right. ,\\
  \tau(E,i)& =\left\{ (e_1,e_2)\ |\ (e_1+ i,e_2)\in E \right\}.
\end{align*}



\begin{figure}
  \begin{center}
  \mbox{\includegraphics{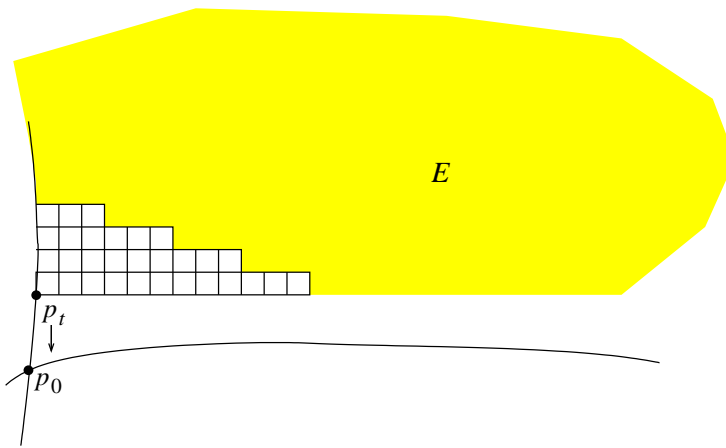}}
    \caption{Example of staircase; the shaded part is $E$, and its
      complement has been drawn as a pile of boxes in staircase
      form. Here $E=\{(e_1, e_2)\in  \Z_{\ge 
        0}^2 | e_1+3e_2\ge 12\}.$ The corresponding ideal is $I_E=(x^3,f)^4$.
      \label{mult4}}
  \end{center}
\end{figure}

The following two propositions are particular cases of the
computations done by \'Evain in \cite{Eva??} and \cite{Eva07}, and we
refer the reader to these works for the proofs.

\begin{Pro}
\label{clm:flatstairs}
  Let $E \subset \Z_{\ge 0}^2$ be a staircase,
  $\m=(m_1,\dots,m_\mu) \in \Z_{> 0}^\mu$ a sequence of positive
  integers, and $I_t=I_E$, $\y=\y_\m$ as defined above. Then 
  \begin{enumerate}
  \item $R_t/I_t$ is flat over $k[[t]]$ and over $k[[y]]$,
  \item $I_t:y_1\dots y_i=I_{\tau(E,\sigma_\m(i))}$,
  \item For every $q\ge p\ge 1$,
  $\HRes{I_t+(t^q)}{x}{p}{}=\HRes{I_t}{x}{p}{}=I_{\tau(E,1)}+(t^p)$, and 
  $\Tr{I_t+(t^q)}{x}{p}{}=\Tr{I_t}{x}{p}{}=(y^{h(E)},x)/(x)$,
\end{enumerate}
\end{Pro}



\begin{Pro}
\label{clm:tracstairs}
  Let $E\subset \Z_{\ge 0}^2$ be a gentle staircase with finite
  complement.
  Then for every couple of integers $q>p\ge 1$,
  \begin{enumerate}
  \item $\tr{I_E+(t^q)}{y}{p}{}=\tr{I_E}{y}{p}{}=h_E(p-1)$,
  \item if $p= \ell_E(i)$ for some $i$ then
  $\Res{I_E+(t^q)}{y}{p}{}=\Res{I_t}{y}{p}{}=(I_{E'}+(t))/(t)$, where
  $E'$ is the only staircase with 
\begin{equation*}
\hat \ell_{E'}(j)=
\begin{cases}
  \hat \ell_{E}(j)-1 & \text {if }j = i\\
  \hat \ell_{E}(j) & \text {if }j \ne i.
\end{cases}
\end{equation*}
  \end{enumerate}
\end{Pro}

\begin{figure}
  \begin{center}
  \mbox{\includegraphics{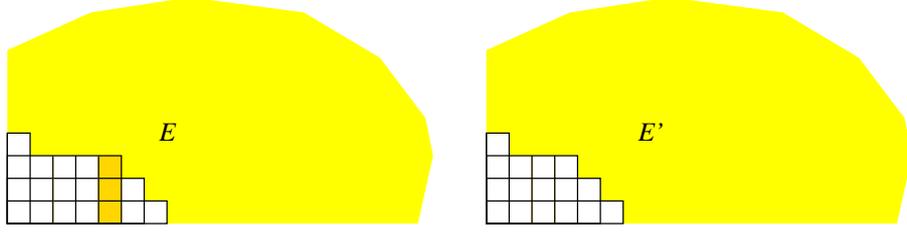}}
    \caption{Example corresponding to proposition
      \ref{clm:tracstairs}; here $p=5$ and $i=2$. The shaded slice in
      the complement to $E$ has to be deleted to obtain $E'$.}
  \end{center}
\end{figure}

$E'$ is the staircase obtained from $E$
by deleting (from its complement) a slice of height
$i+1=h_E(p-1)=\tr{I_E}{y}{p}{}$ and moving everything to the left.

\begin{Cor}
  If $E\subset \Z_{\ge 0}^2$ is a gentle staircase, then
  $\val{I_E+(t^q)}{x}{p}{}=\val{I_E}{x}{p}{}=h(E) $
  and $\val{I_E+(t^q)}{y}{p}{}=\val{I_E}{y}{p}{}=h_E(p-1)$ for
  all $q>p\ge 1$. 
\end{Cor}

\begin{Teo}
\label{clm:stcres}
  Let $E\subset \Z_{\ge 0}^2$ be a given staircase with finite
  complement, and let $\m=(m_1,\dots,m_\mu), 
  \mathbf{tr}=(tr_1, \dots,  tr_\mu)\in \Z_{>0}^{\mu}$ be given
  sequences,
  with $tr_1<tr_2< \dots< tr_\mu$. Define $n_i=\sum_{j<i}(m_j-1)$. Assume that
  \begin{enumerate}
  \item \label{stchat}$\hat \ell_E(tr_i-1) \ge tr_i+1\quad \forall i
    < \mu$,  
  \item \label{stceach} $\ell_E(tr_i-1)-\ell_E(tr_{i+1}-1) \ge h_E(n_i) \quad \forall i<\mu$,
  \item \label{stclast} $\ell_E(tr_\mu-1)> n_{\mu},$
  \item \label{stclasthat} if $\ell_E(tr_\mu)>n_{\mu}$ then $\hat
    \ell_E(tr_\mu-1) \ge tr_\mu+1$.
  \end{enumerate} 
Then there exists $\p=(p_1,\dots, p_m)$ with $m=\sum m_i$ such that
\begin{enumerate}
\item \label{stcres1}for $j=1, \dots, \mu$, $\tr{I_E}{\y_\m}{\p}{j+n_j}=tr_j$,
\item \label{stcres2}for $j+n_j<i\le j+n_{j+1}$,
  $\tr{I_E}{\y_\m}{\p}{i}=h_E(i-j-1)$, 
\item \label{stcres3}$\Res{I_E}{\y}{\p}{m}=(I_{E'})_0$, where
  $E'=\tau(E^{\flat},m-\mu)$ and $E^\flat$ is the staircase with finite
  complement that has
  \begin{align*}
    \hat \ell_{E^{\flat}} (tr_i-1)&=\hat \ell_{E} (tr_i-1)-1, \ i=1, \dots
    \mu,\\
   \hat \ell_{E^{\flat}} (j)&=\hat \ell_{E} (j), \ \text{whenever }j+1
   \not\in \mathbf{tr} .
  \end{align*}
\end{enumerate}
In particular the number of conditions is preserved.
\end{Teo}

$E'$ is the staircase obtained from $E$
by deleting  the leftmost $m-\mu$
slices, and further $\mu$ slices of heights $tr_1, tr_2, \dots,
tr_\mu$. 

\begin{proof} For simplicity  denote $\y=\y_\m$.
  Define $\p$ as follows. 
  $p_{n_j+j}=\ell_E(tr_{j}-1)-n_j$ for $1\le j \le \mu-1$,
  $p_{n_\mu+\mu}=\ell_E(tr_{\mu}-1)-n_\mu$ if $\ell_E(tr_\mu)>n_{\mu}$,
  $p_{n_\mu+\mu}=1$ otherwise;
  $p_{n_j+j+1}=\ell_E(tr_{j}-1)-n_j-h_E(n_j)$ for $1\le j \le \mu-1$
  (so for instance  
  $p_1=\ell_E(tr_1-1)$ and $p_2=\ell_E(tr_1-1)-h(E)$), $p_i=1$ for 
  $i>n_\mu+\mu$ and $p_i=p_{i-1}-1$ for all other $i$. The numerical
  hypotheses \ref{stceach} and \ref{stclast} on the lengths of the
  stairs of $E$ guarantee that with 
  this definition $p_1 > \dots >p_{n_j+j}\ge p_{n_j + j +1} +
  h_E(n_j)\ge \dots \ge p_m$. 
  Then we claim that \ref{stcres1},
  \ref{stcres2} and \ref{stcres3} hold.

\begin{figure}
  \begin{center}
    \mbox{\includegraphics{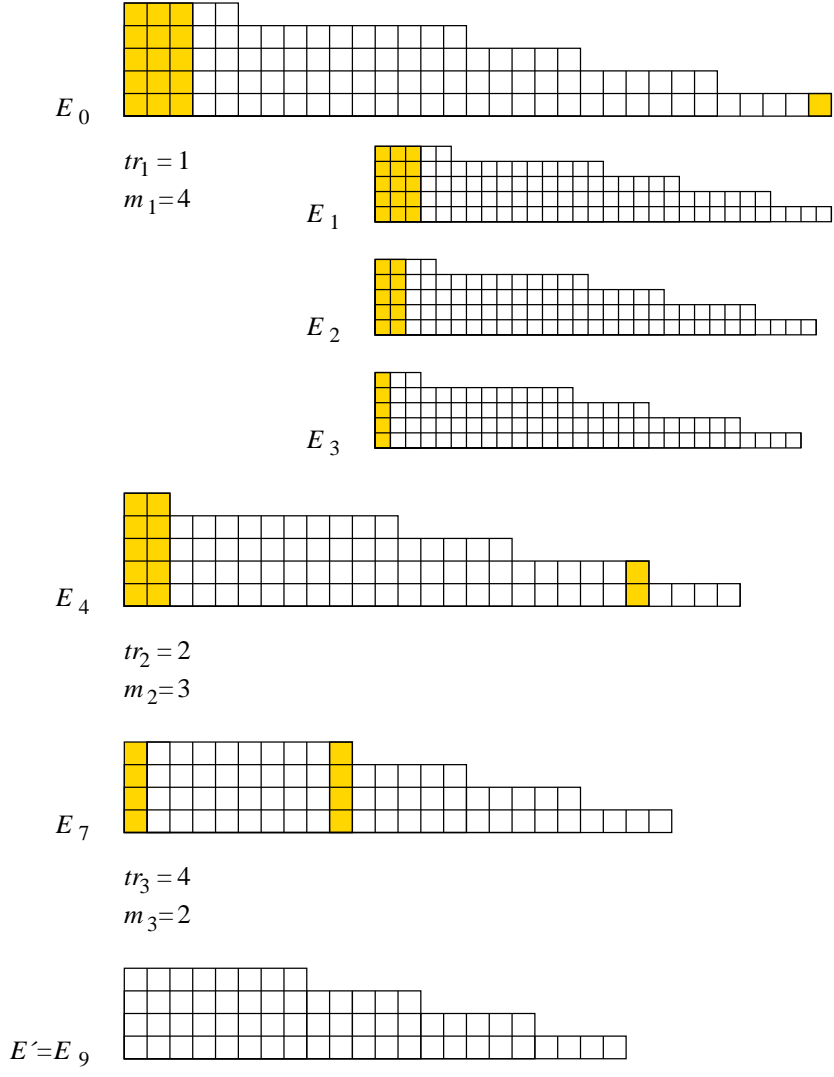}}
    \caption{Example of the computation of residuals in the proof of 
      theorem \ref{clm:stcres}. Given $E=E_0$, $\mathbf{tr}=(1,2,4)$
      and $\m=(4,3,2)$, the complements to the staircases $E_i$ 
      are as shown. For simplicity, starting with $i=4$ we only show 
      the steps $j+n_j$ of the sequence. Shaded, slices to erase.
      \label{sevenst}}
  \end{center}
\end{figure}

To begin with, let us prove claims \ref{stcres1} and 
  \ref{stcres2} for $1\le i< n_\mu+\mu$, and for $i=n_\mu+\mu$ if
  $\ell_E(tr_\mu)>n_{\mu}$. 
Due to proposition \ref{clm:flatstairs}, for
all $i$ and $j$ with $j+n_j<i\le j+1+n_{j+1}$, 
\begin{multline}
\label{eq:stctr}
\val{I_t:(y_1 \dots y_{j+n_j})}{\y_{j+n_j}}{\p_{j+n_j}}{i-j-n_j}
=\val{I_{\tau(E,\sigma_\m(j+n_j))}}{\y_{j+n_j}}{\p_{j+n_j}}{i-j-n_j}=\\
=\val{I_{\tau(E,\sigma_\m(i-1))}}{y_i}{\p_{i-1}}{1}.
\end{multline}
If $i\le j+n_{j+1}$ then (\ref{eq:stctr}) can be evaluated using
proposition \ref{clm:flatstairs}, which gives
\begin{multline*}
 \val{I_{\tau(E,\sigma_\m(i-1))}}{y_i}{\p_{i-1}}{1}=
\tr{I_{\tau(E,i-j-1)}}{y_i}{p_i}{}=\\
=h(\tau(E,i-j-1))=h_E(i-j-1). 
\end{multline*}
On the other hand, if $i=j+1+n_{j+1}$, then using \ref{clm:tracstairs}
we get
\begin{multline*}
\val{I_{\tau(E,\sigma_\m(i-1))}}{y_i}{\p_{i-1}}{1}=\tr{I_{\tau(E,i-j-1)}}{y_i}{p_i}{}=h_{\tau(E,n_{j+1})}(p_{j+1+n_{j+1}}-1)=\\
=h_E(\ell_E(tr_{j+1}-1)-n_{j+1}+n_{j+1}-1)=tr_{j+1}.
\end{multline*}
In both cases the result is bounded above by 
$h_E(n_{j})$  and therefore
$$
\M^{h_E(n_{j-1})}
\subset \Tr{I_t:(y_1 \dots
  y_{j+n_j})}{\y_{j+n_j}}{\p_{j+n_j}}{i-j-n_j}
\subset \Tr{I_t}{\y}{\p}{i},
$$
But the definition of $\p$ gives that 
$p_{j+n_j}-p_i\ge h_E(n_j)+i-j-n_j-1$, and we also have 
$(i-j-n_j)=\sum_{\ell=j+n_j}^{i-1} v(y_\ell)$, therefore
$$
\M^{p_{j+n_j}-p_i+1-\sum_{\ell=j+n_j}^{i-1} v(y_\ell)}
\subset \M^{h_E(n_{j-1})}
$$
and by \ref{clm:cutideal1}, 
$\tr{I_t}{\y}{\p}{i}=\tr{I_t:(y_1 \dots
  y_{j+n_j})}{\y_{j+n_j}}{\p_{j+n_j}}{i-j-n_j}$ is as claimed.

Before considering the cases with $i>\mu+n_\mu$ let us compute
$\Res{I_E}{\y}{\p}{i}$ for $1\le i\le n_\mu+\mu$. We claim that 
 $\Res{I_E}{\y}{\p}{i}=(I_{E_i})_0$, where $E_0=E$,
$E_{j+n_j}$ is the only staircase with finite colength and
\begin{equation*}
\hat \ell_{E_{j+n_j}}(k)=
\begin{cases}
  \hat \ell_{E_{j+n_j-1}}(k)-1 & \text {if }k= tr_j-1\\
  \hat \ell_{E_{j+n_j-1}}(k) & \text {if }k \ne tr_j-1.
\end{cases}
\end{equation*}
and $E_i =\tau(E_{i-1},1)$ for $j+n_j<i\le j+n_{j+1}$. 
Remark that  with this definition, $E_m=E'$ (see figure \ref{sevenst}).
For $i=0$ there is nothing to prove. For $i>0$ we proceed recursively, 
so assume we have proved that $\Res{I_E}{\y}{\p}{i-1}=(I_{E_{i-1}})_0$. 
Now from proposition \ref{clm:cutideal1}, part \ref{ordhres} we deduce that
$\Res{I_t}{\y}{\p}{i}+\M^{p_i-1}=(I_t:(y_1 \dots y_{i-1}))_0+\M^{p_i-1}
=(I_{E_{i-1}})_0+\M^{p_i-1}$, and one has 
$\Res{I_t}{\y}{\p}{i}\subset \Res{I_t}{\y}{\p}{i-1}:y_i$ as well. Therefore 
$\Res{I_t}{\y}{\p}{i}\subset (I_{E_{i-1}})_0+(\M^{p_i-1} \cap (I_{\tau(E_{i-1},1)})_0)=
(I_{E_i})_0$,
where in the case $i=j+n_j$ we use the hypothesis \ref{stchat} (or
hypothesis \ref{stclasthat}) as in the proof of \ref{clm:tracstairs}.
On the other hand, from the first two claims and the key lemma
\ref{clm:cutideal2} we get that $\res{I_t}{\y}{\p}{i}= \# (\Z_{\ge
  0}^2 \setminus E_{i-1}) - 
\tr{I_t}{\y}{\p}{i}= \# (\Z_{\ge 0}^2 \setminus E_{i}) = \dim_k
(R/(I_{E'})_0)$, and we are done. 

Now, for $i>\mu+n_\mu$ and for $i=\mu+n_\mu$ if
$\ell_E(tr_\mu)\le n_{\mu}$, we have $p_i=1$, hence by \ref{clm:flatstairs} 
$$\HRes{I_E}{\y}{\p}{i}=\Res{I_E}{\y}{\p}{i} k[[t]]+(t)
=I_{\tau(E_{\mu+n_\mu},i-\mu+n_\mu)}+(t)$$ so
$\Res{I_E}{\y}{\p}{i}=(I_{E_i})_0$ in these cases as well (in
particular the last claim follows) and
$\tr{I_E}{\y_m}{\p}{i}=h_E(i-\mu-1)$ as claimed.
\end{proof}

\section{Proof of theorem \ref{clm:regsmalle}}
\label{sec:applications}


The sequence of specializations to which we apply the preceding technique
in order to prove theorem \ref{clm:regsmalle} was already introduced 
in \cite{Roe01} and used in
\cite{Roe01b}, \cite{Roe01a} and \cite{Roe06}. It consists in
introducing satellite points: one
first specializes each point to be infinitely near to the previous
one and then, step by step, the third point is brought to the
first irreducible exceptional component (at its intersection point
with the second), then the fourth, and so on. 
As a byproduct, the result we obtain is slightly stronger,
as it shows regularity of linear systems defined by a more general
class of \emph{cluster} schemes $Z$.
Their ideals are obtained as follows. If
$\P^2=X_0\leftarrow X_1 \leftarrow X_2 \leftarrow \dots \leftarrow X_n$
is a sequence of blowing ups centered at the points $p_1, \dots, p_n$,
$\pi$ denotes their composition, and $E_i$ is the exceptional divisor
in $X_n$ above $p_i$ for each $i$, then  
$$\mathcal{I}=\pi_*(\O_{X_n}(-e E_1-\dots -e E_n))$$
is the ideal sheaf defining the cluster scheme that has each point
$p_i$ with multiplcity $e$. We assume the clusters are consistent, i.e., 
$\O_{X_n}(-e E_1-\dots -e E_n)$ cuts non-negatively each
irreducible component of the $E_i$. Under these hypotheses there can
be no satellites among the $p_i$, i.e., each point
belongs to at most one irreducible exceptional component.
The information of proximities satisfied by the points (i.e., which
points belong to which exceptional components) can be encoded into
Enriques diagrams or proximity matrices \cite{Cas00}. Clusters with the same
diagram $\D$ (or with the same matrix) are parameterized by an
irreducible quasiprojective variety $\Cl(\D)$  \cite{Roe04b}, and
the expression ``general clusters with diagram $\D$'' refers to
clusters parameterized by a Zariski open subset of $\Cl(\D)$. 
\begin{Teo}
\label{clm:smallmult}
 Let $n,e$ be positive integers, with $\sqrt{n} \ge 2e$ and $e>2$. Then for every 
 consistent weighted Enriques diagram $\D$ with exactly $n$
 vertices, all of multiplicity $e$, general
 weighted clusters of type $\D$ on the projective plane have
 maximal rank in all degreees. 
\end{Teo}
The technical hypothesis $e > 2$ is not really restrictive: for
$e\le 2$ the result is known to be true with no restriction on $n$ (see
\cite{Roe06}). Theorem \ref{clm:regsmalle} corresponds to the particular 
in which $\D$ consists of $n$ distinct points.

Section \ref{sec:gentle} introduces the specializations that will be
used to prove theorem \ref{clm:smallmult}. These consist in
gradually increasing the contact of suitable schemes defined by
monomial ideals with a curve of selfintersection $-1$, which in the
application will be the exceptional divisor of blowing up a
point. Using some additional blow ups and the results of section
\ref{sec:method} we show how to bound the desired limits. Then we
exploit results from \cite{Roe06} to prove the theorem in section 
\ref{sec:manypoints}. 

\subsection{Gentle staircases on blown-up surfaces}
\label{sec:gentle}

Let $S$ be a smooth, projective, algebraic surface, and $D$ a
$(-1)$-curve on it. 
Given $p\in D$, let $R$ be the completion of the local ring $\O_{S,p}$,
and fix an isomorphism $R=\hat \O_{S,p}\cong k[[x,y]]$
such that $y=0$ is a local equation for $D$ at $p$.
If no confusion is likely, the maximal ideals of $\O_{S,p}$
and $k[[x,y]]$ will be both denoted by $\M$.
By the second part of lemma \ref{clm:indepx}, every possible monomial
ideal $I_{E,f,g}\subset k[[x,y]]$
can be defined by 
$f, g\in \O_{S,p}$. $I_{E,f,g}\cap\O_{S,p}$ is primary with respect to
the maximal ideal too, of the same colength $\#(\Z_{\ge 0}^2 \setminus E)$. 

Let now  $C\subset S$ a curve through $p$, and $g\in  \O_{S,p}$ a
local equation for $C$. If $E$ is a gentle staircase, 
$I_{E,f,g}$ does not depend on $f$ by \ref{clm:indepx}; then we denote
$I_{E,g}=I_{E,f,g}$ and define $\mathcal{I}_{p,E,C}$ or
$\mathcal{I}_{p,E,g}$ to be the ideal sheaf with
cosupport at $p$ and stalk $I_{E,g}\cap\O_{S,p}$, and
$Z_{p,E,C}\subset S$ or $Z_{p,E,g}\subset S$ the zeroscheme it
defines.


%
If $L$ is a divisor such that $\len (Z_{p,E,C}\cap D) > L \cdot D$ then
$D$ is a fixed part of all 
curves in $|L|$ that contain $Z_{p,E,C}$, if they exist. 
In order to compute $\len (Z_{p,E,C}\cap D)$
we introduce a couple of definitions.
 Given an integer  $r\ge 0$, we denote 
$h^r_E=\min\{i|\hat \ell_E(i)\le r\}$. 
 If the function $\hat \ell_E:\Z_{\ge 0}
 \longrightarrow \Z_{\ge 0}$ is  non-increasing in the interval
 $[h^r_E, \infty)$ then we say that $E$ is $r$-gentle. Observe that
for $r>0$ $r$-gentle implies $(r-1)$-gentle, and $0$-gentle implies
gentle.

\begin{Lem}
  Let  $\pi:S'\longrightarrow S$ be the blow-up of $p$, let $D_p$ be
the exceptional divisor and $\tilde D$ the strict transform of
$D$. Let $C$ be a curve going through $p$ and smooth at $p$, $\tilde
C$ its strict transform on the blow-up and $p'=\tilde C \cap D_p$.
If $E$ is a 1-gentle staircase with finite complement, then
  $\mathcal{I}_{p,E,C}=\pi_*(\mathcal{I}_{p',\tilde E,\tilde C}
  \otimes \mathcal{O}_{S'}(-h(E) D_p) )$, 
where $\tilde E$ is the staircase with finite complement that has
$$\hat\ell_{\tilde E}(i)=\max\{\hat \ell_E(i)-1,0\}$$
for all $i$. Moreover, if $E$ is $r$-gentle $r\ge 1$, then $\tilde E$
is $(r-1)$-gentle. 
\end{Lem}
\begin{proof}
  Both ideal sheaves $\mathcal{I}_{p,E,C}$ and
  $\pi_*(\mathcal{I}_{p',\tilde E,\tilde C}
  \otimes \mathcal{O}_{S'}(-h(E) D_p))$ have cosupport at $p$, and
  their stalk there is 
  primary with respect to the maximal ideal. Therefore, it will be
  enough to see that their extensions to the completion of $\O_{S,p}$
  coincide. 

  Let $f=0$ is a local equation of
  $C$, $y=0$ a local equation of $D$ and assume that  the isomorphism
  $\hat \O_{S,p}\cong k[[x,y]]$ has been chosen in such a way that
  $(x,f)=\M$. Then $\hat \O_{S,p}\cong k[[x,f]]$, and $\hat
  \O_{S',p'}\cong R' = R[[f/x]]=  k[[x,f/x]]$. Therefore $f/x=0$ is a local
  equation of $\tilde C$, $x=0$ 
  a local equation for $D_p$ and, if $\tilde D$ goes through $p'$
  (which means that $(D\cdot C)_p>1$) then
  $y/x$ is a local equation of $\tilde D$.
  Then, the stalk at $p$ of $\pi_*(\mathcal{I}_{\tilde E,\tilde C}
  \otimes \mathcal{O}_{S'}(-h(E)  D_p)) $ is 
  $\O_{S,p}\cap (x^{h(E)} I_{\tilde E,\tilde C})$, and its
  extension to $R$ is
  \begin{gather*}
   R\cap (x^{h(E)} I_{\tilde E,\tilde C})=
   k[[x,f]]\cap
   \left(\left(x^{h(E)+e_1}(f/x)^{e_2}\right)_{(e_1,e_2)\in \tilde E'}
    \right)= \\ =
   \left(x^{h(E)+e_1-e_2}f^{e_2}\right)_{(e_1,e_2)\in \tilde E'} 
   \subset k[[x,f]]=R. 
  \end{gather*}
  Now, if $E$ is 1-gentle, it is immediate to check that
  $(h(E)+e_1-e_2,e_2) \in E$ if and only if $(e_1,e_2)\in \tilde E'$, and we
  are done. The claim on $(r-1)$-gentleness of $\tilde E'$ is immediate
  from the definitions.
\end{proof}
We can iterate this process, by blowing-up $p'$,  then $p''=\tilde
C \cap D_{p'}$, and so on, if the staircase is gentle enough:

\begin{Cor}
\label{clm:rgentle}
  Let $E$ be a $r$-gentle staircase with finite complement, 
$C$  a curve going through $p$ and smooth at $p$, and $(C \cdot D)_p=s
\ge r$. Let $p, p', \dots, p^{r}$ be the
  first $r+1$ points on $C$ infinitely near to $p$, and $S_r$ the
  surface obtained by blowing up $p, p', \dots, p^{r-1}$, on which
  $p^{r}$ lies. 
  Then   $\mathcal{I}_{p,E,C}=\pi_*(\mathcal{I}_{p^{r},\tilde E,\tilde
  C}\otimes \O_{S_r}(-d_0D_p-d_1D_{p'}-\dots-d_{r-1}D_{p^{r-1}}))$,
  where $\tilde E$ is the staircase with finite complement that has 
$$\hat\ell_{\tilde E}(i)=\max\{\hat \ell_E(i)-r,0\}$$
for all $i$, $d_i=\max\{j|\hat\ell_E(j)>i\}+1$ and $D_{p^i}$ denotes the
total transform in $S_r$ of the $i$th exceptional divisor. If moreover
$\lmin(E)\ge r$ then $d_i=h_E(i)$.
\end{Cor}

\begin{Cor}
\label{clm:indepf}
  Let $E\subset \Z_{\ge 0}^2$ be a staircase with finite
  complement, and $C$  a curve going through $p$ and smooth at $p$
  such that $(C \cdot D)_p=r$.
  If $E$ is $r$-gentle then $\len (Z_{p,E,C}\cap
  D)=\ell_E(h^r_E)+rh^r_E$. If moreover $\lmin(E)\ge r$ then 
  $\len (Z_{p,E,C}\cap D)=\sum_{i=0}^{r-1} h_E(i)$.
\end{Cor}

\begin{proof}
  Let as before $p, p', \dots, p^{r}$ be the
  first $r+1$ points on $C$ infinitely near to $p$, and $S_r$ the
  surface obtained by blowing up $p, p', \dots, p^{r-1}$. $p^{r}\in S_r$
  does not belong to the strict transform of $D$ because $(C \cdot
  D)_p=r$, and due to Corollary 
  \ref{clm:rgentle}, $\mathcal{I}_{p,E,C}=\pi_*(\mathcal{I}_{p^{r}\tilde
  E,\tilde f}\otimes
  \O_{S_r}(-d_0D_p-d_1D_{p'}-\dots-d_{r-1}D_{p^{r-1}}))$, where
  $d_i=\max\{j|\hat\ell_E(j)>i\}+1$. Then by the projection formula,
  $\len (Z_{p,E,C}\cap D)= \sum_{i=0}^{r-1} d_i=\ell_E(h^r_E)+rh^r_E$,
  as wanted. 

  If moreover $\lmin(E)\ge r$ then $d_i=h_E(i)$ so by the projection
  formula again  $\len (Z_{p,E,C}\cap
  D)=\sum_{i=0}^{r-1} h_E(i)$. 
\end{proof}

\begin{Rem}
  \label{clm:rkr+1} It is worth noting that if $\lmin(E)\ge r+1$ then
  $E$ is $r$-gentle and the previous two corollaries apply.
\end{Rem}

  Given a triple $(L, E, r)$, where $L$ is a
 divisor class on $S$, $r$ is a positive integer and $E$ is a
 $r$-gentle staircase with finite complement, 
 we say that a linear system $\Sigma$ on $S$ has type
 $(L, E, r)$ if there is a curve $C$ 
 through $p$, smooth at $p$ and with $(C\cdot D)_p=r$, 
such that $\Sigma=\P(H^0(\mathcal{I}_{p,E,C}\otimes L))$.
If  $L\cdot D \ge \ell_E(h^r_E)+rh^r_E$, then 
the type $(L,E,r)$ is called \emph{consistent}.

Given a family of curves $C_t$ through $p$, 
the intersection number $(C_t\cdot D)_p$ may depend on the parameter
$t$, i.e., one may have 
$(C_t \cdot D)_p=r$ and $(C_0\cdot D)_p=r+1$, for instance. Then one obtains a
family of linear systems $\Sigma_t$, $t\ne 0$ of type $(L, E, r)$
whose limit when $t\mapsto 0$ is of different type.

\begin{Lem}
\label{clm:trivstep}
    Let $(L,E,r)$ be a consistent type with $r>1$. Every linear system of type
  $(L, E, r)$ contains as a sublinear system the moving part of the
  limit of a family of linear systems of type   $(L,E,r-1)$.
\end{Lem}
\begin{proof}

  Fix local coordinates $(x,y)$ such that $y=0$ is a local equation
for $D$. For every curve $C$ on $S$ going through $p$ with $(C\cdot D)_p=r$,
let $f\in \O_{S,p} \subset \hat \O_{S,p} =k[[x,y]]$ be a local equation
for $C$.   

For  $t\ne 0$, $f+tx^{r-1}$ is a local equation at $p$ of a curve
$C_t$ with $(C_t \cdot D)_p=r-1$. Define
$\Sigma_t=\P(H^0(\mathcal{I}_{p,E,C_t}\otimes L))$.
Then it is clear by the definitions that
\begin{equation*}
 \lim_{t\mapsto 0} \Sigma_t \subset
 \Sigma:=\P(H^0(\mathcal{I}_{p,E,C}\otimes L)). \qedhere
\end{equation*}
\end{proof}

\begin{Teo}
\label{clm:steptypes}
  Let $(L,E,r)$ be a consistent type such
  that $\lmin(E)\ge r+h^{r}_{E}+1$ and $h^r_E\ge 2$. There exist an
  integer $\mu\ge 0$ and an
  $(r+1)$-gentle staircase $E'$ with
  \begin{enumerate}
  \item \label{stepcodim}$\#\left(\Z_{\ge 0}^{2}\setminus E'\right)+\mu(L\cdot
D)+\binom{\mu+1}{2}=\#\left(\Z_{\ge 0}^{2}\setminus E\right),$
  \item \label{steptotal} $\tau(E,\mu r)\subset E' \subset \tau(E,\mu(r+1))$,
  \item \label{steplength} $\ell(E')=\ell(E)-\mu(r+1)$, 
  \item \label{stepstairs} $\lmin(E')\ge \lmin(E)-1$,
  \item \label{steptop} if $\ell_E(h(E)-1)> \mu r+1$ and $\mu \ge 1$
    then $h(E')=h(E)$ and $\ell_{E'}(h(E')-1)= \ell_E(h(E)-1)-(\mu r+1)$,
  \end{enumerate}
such that $(L-\mu D, E', r+1)$ is consistent and every linear system of type
  $(L-\mu D, E', r+1)$ contains as a sublinear system the moving part
  of a limit of linear systems of type 
  $(L,E,r)$. 
\end{Teo}

\begin{proof}
For every integer $i>0$, consider the following
quantities:
\begin{align*}
  s_i&=\sum_{j=r(i-1)}^{ir-1}h_E(j), \\
  tr_i&=L\cdot D + i -s_i, \\
\end{align*}
It is clear that $s_1 \ge s_2 \ge \dots$, and therefore
$tr_1<tr_2<\dots$. Let $\mu=\max\{i\ | \ \ell_E(tr_i -1)>ri\}$,
$\mathbf{tr}=(tr_1, tr_2, \dots, tr_\mu)$. If $\mu=0$ then either
$tr_1>h(E)$ or $tr_1=h(E)$ and $\ell_E(h(E)-1)<r$, in which case
$h^r_E=h(E)-1$; in both cases $(L,E,r+1)$ is consistent and the claims
follow from lemma \ref{clm:trivstep} setting $E'=E$.

So assume $\mu \ge 1$. We claim that the staircase
$E'=\tau(E^{\flat},\mu r)$, satisfies the stated conditions, where
$E^\flat$ is the staircase with finite complement that has
  \begin{align*}
    \hat \ell_{E^{\flat}} (tr_i-1)&=\hat \ell_{E} (tr_i-1)-1, \ i=1, \dots
    \mu,\\
   \hat \ell_{E^{\flat}} (j)&=\hat \ell_{E} (j), \ \text{whenever }j+1
   \not\in \mathbf{tr}.
  \end{align*}
$E'$ is obtained from $E$
by deleting the ``leftmost'' $\mu r$
slices, and further $\mu$ slices of heights $tr_1, tr_2, \dots,
tr_\mu$ (the hypothesis on $\lmin$ guarantees that such slices exist
and that this description is correct). Thus, claims \ref{stepcodim},
\ref{steptotal}, \ref{steplength}, \ref{stepstairs} and \ref{steptop} 
follow.

Moreover, as $\lmin(E')\ge \lmin(E)-1\ge r+h^r_E>r+1$, $E'$ is
$(r+1)$-gentle and for every $C$ through $p$ with $C \cdot D=r+1$,
$\len(Z_{p,E',C} \cap D) = \sum_{i=0}^{r}h_{E'}(i)$, which by
\ref{steptotal} is at most equal to $\sum_{i=\mu r}^{(\mu+1)r}h_{E}(i)$ and by
the definition of $\mu$ this is at most $\L \cdot D +\mu$. Therefore
$(L-\mu D, E', r+1)$ is consistent.

It remains to be seen that every linear system of type
  $(L-\mu D, E', r+1)$ contains as a sublinear system the moving part
  of a limit of linear systems of type $(L,E,r)$ (the fixed part being
  $\mu D$).
So let $C$ be a curve on $S$ going through $p$, with $(C\cdot D)_p=r+1$,
and assumme that local coordinates $(x,y)$ have been chosen in
$\O_{S,p}$ in such a way that $y=0$ is a local equation for $D$, and
$f\in \O_{S,p} \subset \hat \O_{S,p} =k[[x,y]]$ is a local equation
for $C$. We need to prove that
$\Sigma=\P(H^0(\mathcal{I}_{p,E',C}\otimes (L-\mu D))$ 
contains as a sublinear system the moving part of a limit of linear
systems of type $(L,E,r)$.  

For every $t \ne 0$, $f+tx^r$ is a local equation at $p$ of a curve
$C_t$ with $(C_t \cdot D)_p=r$. Define
$\Sigma_t=\P(H^0(\mathcal{I}_{p,E,C_t}\otimes L))$.
We claim that
\begin{equation}
  \label{eq:desired}
 \lim_{t\mapsto 0} \Sigma_t \subset \Sigma+ \mu D. 
\end{equation}
The first $r$ points on $C_t$ infinitely
near to $p$ lie on $D$ as well, so they do not depend on $t$; denote them
$p, p', \dots, p^{r-1}$, and let $\pi:S_r \rightarrow S$ 
be the blowing up of these points. The $(r+1)$th point on
$C_t$ infinitely near to $p$ depends on $t$; let it be $p_t^r\in
S_r$. We shall compute the limit of the $\Sigma_t$ on $S_r$ rather
than on $S$. Indeed, corollary \ref{clm:rgentle} shows
that $\mathcal{I}_{p,E,C_t}=\pi_*(\mathcal{I}_{p_{t}^{r},\tilde E,\tilde
  C_t}\otimes
\O_{S_r}(-d_0D_p-d_1D_{p'}-\dots-d_{r-1}D_{p^{r-1}}))$, where $\tilde
E$ is the staircase with finite complement that has 
$$\hat\ell_{\tilde E}(i)=\max\{\hat \ell_E(i)-r,0\}$$
for all $i$ and $d_i=h_E(i)$. Setting
$\Sigma'_t=\P(\mathcal{I}_{p_{t}^{r},\tilde E,\tilde C_t}\otimes
\O_{S_r}(L-d_0D_p-d_1D_{p'}-\dots-d_{r-1}D_{p^{r-1}}))$, it is clear
that $\Sigma_t \overset{\pi_*}\cong \Sigma'_t$, and 
\begin{equation}
  \label{eq:doitonblowup}
\lim_{t\mapsto 0} \Sigma_t=\pi_*\left(\lim_{t\mapsto 0} \Sigma'_t\right).  
\end{equation}

On $S_r$, the point $p^r_0$ belongs to the strict transforms $\tilde
C$ and $\tilde D$ of the curves $C$ and $D$ respectively, and to the
exceptional divisor $D_{p^{r-1}}$; at $p^r_0$, $\tilde D$,
$D_{p^{r-1}}$ and $\tilde C$ are pairwise transverse. Thus, there
exist $x_r, y_r \in \hat\O_{S_r,p^r_0}$ local parameters such that
$y_r=0$, $x_r=0$, and $x_r+y_r=0$ are local equations of $\tilde D$,
$D_{p^{r-1}}$ and $\tilde C$ respectively. Then $f_t=x_r+y_r+t=0$, for $t$
in a neighbourhood of 0, is an equation of $C_t$ in a neighbourhood of $p^r_0$.

  Let $I_{\tilde E}=\left(x_r^{e_1}f_t^{e_2}\right)_
{(e_1,e_2)\in \tilde E},$
as in section \ref{sec:monomial}. Define also for $i=1, \dots, \mu$,
\begin{align*}
  l_i&=\ell_E(h_E(ri)-1),\\
  m_i&=\min\{r,ri-l_{i-1},l_i -ri\},
\end{align*}
$\mathbf{m}=(m_1, m_2, \dots, m_\mu)$, $m=\sum m_i$, and
$n_i=\sum_{j<i}(m_j-1)$.


\begin{figure}
  \begin{center}
    \mbox{\includegraphics{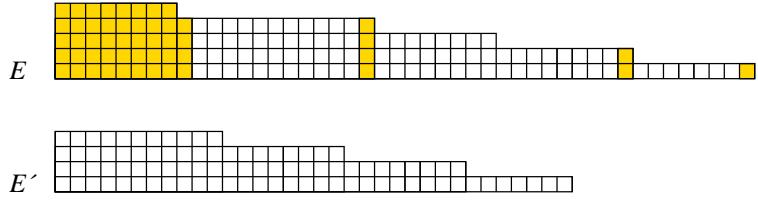}}
    \caption{Example of the computation of $E'$ in 
      theorem \ref{clm:steptypes}. The figure shows the complement to
      a given staircase $E$; the shaded part has to be erased to
      obtain the complement to $E'$ if $r=3$ and $L\cdot D=15$. 
      Note that the staircases $\tilde E$ and $\tilde E'$ of
      the proof coincide, in this example, with the staircases $E_0$
      and $E_9$ shown in figure \ref{sevenst}.
      \label{onest}}
  \end{center}
\end{figure}

It is not difficult to check that the hypotheses of theorem \ref{clm:stcres}
are satisfied for ideals defined by the
staircase $\tilde E$, with
$\y=\y_{\m}=(y,x,\dots,x,y,x,\dots,x,y,x,\dots)$ as in section
\ref{sec:monomial}. Moreover, the staircase
$\tilde E'$ given by \ref{clm:stcres} satisfies
$\hat\ell_{\tilde E'}(i)=\max\{\hat \ell_{E'}(i)-r,0\}$.
Therefore by \ref{clm:rgentle} one obtains
$\mathcal{I}_{p,E',C}=\pi_*(\mathcal{I}_{p_{0}^{r},\tilde E',\tilde 
  {C}}\otimes
\O_{S_r}(-d'_0D_p-d'_1D_{p'}-\dots-d'_{r-1}D_{p^{r-1}}))$ for
  $d'_i=h_{E'}(i)$. In
particular 
\begin{multline}
  \label{eq:blowuptrans}
\P(H^0(\mathcal{I}_{p,E',C}\otimes \O_{S}(L-\mu D)))= \\ =
\P(\pi_*H^0(\mathcal{I}_{p_{0}^{r},\tilde E',\tilde {C}}\otimes
\O_{S_r}(L-\mu D-d'_0D_p-d'_1D_{p'}-\dots-d'_{r-1}D_{p^{r-1}})) . 
\end{multline}

Let $V\subset \hat \O_{S_r,p^r_0}$ be the
image of the natural morphism
$H^0(\O_{S_r}(L-d_0D_p-d_1D_{p'}-\dots-d_{r-1}D_{p^{r-1}}))\rightarrow
\hat \O_{S_r,p^r_0}$. Lemma \ref{clm:resclasses} below shows that 
\begin{equation}
  \label{eq:enoughtrace}
 \frac{\Res{V}{y_1 \cdots y_{i-1}}{}{}}{\Res{V}{y_1 \cdots
    y_{i-1}}{}{} \cap (y_i)} \longrightarrow
\frac{R/(y_i)}{\Tr{I_t}{\y}{\p}{i}} 
\end{equation}
is injective for $i=1, \dots, \mu$, so theorem
\ref{clm:formal} applies as well, and therefore 
\begin{multline}
  \label{eq:limit}
 \lim_{t\mapsto 0} \Sigma'_t= \mu D + (d'_0-d_0) D_p + \dots
 (d'_{r-1}-d_{r-1}) D_{p^{r-1}} + \\
+\P(H^0(\mathcal{I}_{p_0^r,\tilde
  E',\tilde f}(L-\mu
D-d'_0D_p-d'_1D_{p'}-\dots-d'_{r-1}D_{p^{r-1}}))). 
\end{multline}
Now it suffices to put
(\ref{eq:doitonblowup}), 
(\ref{eq:blowuptrans}) and  (\ref{eq:limit}) together to see that
  (\ref{eq:desired}) holds.
\end{proof}

\begin{Lem}
  \label{clm:resclasses} For
each $i$ there is a divisor class $F_i$ on $S_r$ such that
\begin{enumerate}
\item \label{resclass} $\Res{V}{y_1 \cdots y_{i-1}}{}{}\subset \hat \O_{S_r,p^r_0}$ is the
  image of the natural  morphism 
 $$H^0(\O_{S_r}(F_i))\overset{\rho_i}\longrightarrow \hat \O_{S_r,p^r_0},$$
\item \label{restrace} $F_i \cdot E_i<\tr{I_t}{\y}{\p}{i}$, where $E_i$ is the irreducible divisor
  defined locally by $y_i=0$, that is, $E_i=\tilde D$ if $i=n_j+j$ for some
  $j$, and $E_i=D_{p^{r-1}}$ otherwise.
\end{enumerate}
\end{Lem}
\begin{proof}
Let us define the $F_i$ by recurrence on $i$. To begin
with, set $F_1=L-d_0D_p-d_1D_{p'}-\dots-d_{r-1}D_{p^{r-1}}$. By
definition and assuming $F_{i-1}$ satisfies the claims,
it is clear that $\Res{V}{y_1 \cdots y_{i-1}}{}{}$ is the 
image of the natural morphism
$$H^0(\O_{S_r}(F_{i-1}-E_{i-1}))\longrightarrow \hat \O_{S_r,p^r_0}.$$
However, the divisor class $F_{i-1}-E_{i-1}$ need not be consistent, 
i.e., it may intersect negatively some irreducible components of the divisors
$D_{p^j}$, which in that case become fixed parts of $|F_{i-1}-E_{i-1}|$. We
define $F_{i}$ to be the 
consistent system obtained from $F_{i-1}-E_{i-1}$ by unloading (i.e.,
subtracting the fixed divisors). Clearly $H^0(\O_{S_r}(F_{i}))\cong
H^0(\O_{S_r}(F_{i-1}-E_{i-1}))$, the isomorphism being given by the
subtraction of the fixed divisors, which do not pass through $p^r_0$,
and therefore  $\Res{V}{y_1 \cdots y_{i-1}}{}{}\subset \hat
\O_{S_r,p^r_0}$ is the image of $\rho_i$ as stated.

Now compute $F_i$ and $F_i
\cdot E_i$. Let $j$ be such that $n_j+j+1 \le i \le n_{j+1}+j+1$, and
define 
\begin{equation*}
k_0 =
\begin{cases}
  0 & \text{if }h_E(r(j-1))=h_E(rj),\\
\ell_E(h_E(r(j-1))-1)-r(j-1)& \text{if }h_E(r(j-1))>h_E(rj),
\end{cases}
\end{equation*}
and for $k=0, \dots, r-1,$
 \begin{equation*}
  d_k^i=
\begin{cases}
    h_E(jr) & \text{if }k-k_0 \le i-(n_j+j+1), \\
    h_E(jr)-1 & \text{if }k-k_0> i-(n_j+j+1).
  \end{cases}
\end{equation*}
Remark that the definitions of $d_k^i$ and $m_i$ imply that, if $i=n_j+j$ then
$d_k^i=h_E(r(j-1)+k)$, and if $i=n_j+j+1$ then $d_k^i=h_E(r(j-1)+k)-1$. 
Now it is an elementary unloading exercise to show that
  \begin{equation}
    \label{eq:descfi}
F_i=L - jD - d_0^iD_p - \dots -d_{r-1}^i
D_{p^{r-1}}.
  \end{equation}

Finally, observe that by (\ref{eq:descfi}), if $i=n_{j}+j$ then $F_i
\cdot E_i= F_i \cdot \tilde D=(L\cdot D) + j$ and for all other values
of $i$, $F_i \cdot E_i= F_i \cdot D_{p^{r-1}}=h_E(jr)-1$, whereas if $i=n_{j}+j$
then $\tr{I_t}{\y}{\p}{i}=(L\cdot D) + j+1$ and for all other values
of $i$,  $\tr{I_t}{\y}{\p}{i}=h_E(i-j)>h_E(jr)-1$, so $F_i \cdot
E_i<\tr{I_t}{\y}{\p}{i}$.
\end{proof}

\subsection{Equimultiple clusters with many points}
\label{sec:manypoints}

\begin{figure}
  \begin{center}
    \mbox{\includegraphics{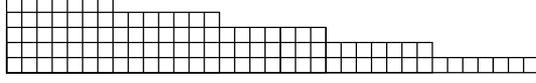}}
    \caption{The complement to staircase $E_1(8,5)$.
      \label{e85}}
  \end{center}
\end{figure}

In this section, $\pi:S\rightarrow\P^2$ is the blowing up of $\P^2$ at
a point, $D$ is the exceptional divisor and 
$(n,e)$ is a couple of integers with $2<e\le\sqrt{n}/2$.

Define $E_1(n,e):=\{(e_1,e_2) \in \Z^{2}_{\ge 0}| 
e e_1+(n-1)e_2\ge e(n-1)\}$. In other words, $E_1(n,e)$ is the staircase
of height $e$ and $\hat \ell_{E_1(n,e)}(i)=n-1$ for
$0\le i < e$. Let $L$ be a divisor class with
$L\cdot D=e$, such as $\pi^*(dH)-eD$, where $H$ is the class of a line
in the plane and $d$ a positive integer. Then $(L, E_1(n,e),1)$ is a
consistent type satisfying 
the requirement of theorem \ref{clm:steptypes}; let $\mu_1$ and
$E_2(n,e)$ be the integer and staircase given by \ref{clm:steptypes}. 
$(L-\mu_1 D, E_2(n,e),2)$ is again a consistent type satisfying
the requirement of theorem \ref{clm:steptypes}; let $\mu_2$ and
$E_3(n,e)$ be the corresponding integer and staircase. As long as the
hypotheses of the theorem are satisfied, we keep using it to define
integers $\mu_3, \mu_4, \dots$ and staircases $E_4(n,e), E_5(n,e),
\dots$ When no confusion may arise, we denote the staircases simply
$E_1, E_2, \dots$ and we also denote $L_2=L-\mu_1D$, $L_3=L_2-\mu_2
D$, \dots Let $\lastr(n,e)$ be the last $r$ such that $E_r$ is
defined, i.e., $(L_{\lastr(n,e)},E_{\lastr(n,e)}(n,e),\lastr(n,e))$ is a
consistent type and either 
$h^{\lastr(n,e)}_{E_{\lastr(n,e)}(n,e)}\le 1$ or 
$\lmin(E_{\lastr(n,e)})\le \lastr(n,e)+h^{\lastr(n,e)}_{E_{\lastr(n,e)}(n,e)}$.

\begin{Lem}
\label{clm:rmin}
  Let $r$ be a positive integer, and denote
  $M=\sum\mu_i$, with the summation running over all $i\le \min\{r-1,
  \lastr(n,e)\}$.
  \begin{enumerate}
\item \label{erdefined}
if $r\le \lastr(n,e)$, then 
  \begin{enumerate}
  \item \label{erdefined-hatell}$\lmin(E_r)\ge n-r$,
  \item \label{erdefined-r}$(r-1)(r+2)e \ge 2(n-1)(e-h(E_r))$, and
  \item \label{erdefined-codim}$\#(E_r \setminus
 E_1)=eM+\binom{M+1}{2}$,
  \end{enumerate}
  \item \label{smallr} If $\binom{r}{2}e+r < n-1$, then $r\le \lastr(n,e)$,
  and $\mu_1=\mu_2=\dots=\mu_{r-1}=e$,
 \item \label{bigrcodim} if $\binom{r}{2}e+r \ge n-1$ and $r\le \lastr(n,e)$,
  then $\#(E_r \setminus E_1)\ge e(n-1)/2$.
  \end{enumerate}
\end{Lem}

\begin{proof}
Observe first that due to \ref{clm:steptypes}, claim \ref{stepstairs}, for
every $1<r\le \lastr(n,e)$ one has $\lmin(E_r)\ge\lmin(E_{r-1})-1$, and
therefore $\lmin(E_r)\ge n-r$.
Now because of \ref{clm:steptypes}, claim \ref{steptotal}, $E_r\supset
  \tau(E_1,(e-h(E_r))(n-1))$. We also have $E_r \subset \tau(E_1,\sum_{i=1}^{r-1}
  \mu_i(i+1))$ and $\mu_i \le e$, hence
$$\frac{(r-1)(r+2)}{2}e \ge (n-1)(e-h(E_r)).$$ 
Thus, claim \ref{erdefined} follows from because of
\ref{clm:steptypes}, claim \ref{stepcodim}: 
$$\#(E_r \setminus
 E_1)=\binom{e+M+1}{2}-\binom{e+1}{2}=eM+\binom{M+1}{2}.$$
%
Now because claim \ref{erdefined-hatell} holds, 
$r+1\le \lastr(n,e)$ whenever $2r \le n-e-1$, and in
particular $r\le \lastr(n,e)$ whenever $\binom{r}{2}e+r < n-1$; in
such a case moreover, due to 
\ref{clm:steptypes}, claim \ref{steptop}, 
$\mu_1=\mu_2=\dots=\mu_{r-1}=e$. 

Finally for the third claim, let $r_0$ be the maximal integer with
$\binom{r_0}{2}e+r_0 < n-1$. The hypothesis says $r\ge r_0$, and
then $M=\mu_1+\mu_2+\dots+\mu_{r-1}\ge e(r_0-1)$. 
Then by claim \ref{erdefined-codim},
$$\#(E_r \setminus E_1)=eM+\binom{M+1}{2},$$
which is not less than $e(n-1)/2$ by the definition of $r_0$ and the
inequality $M\ge (r_0-1)e$.
\end{proof}

\begin{Lem}
\label{clm:goesto2}
For every $r\le \min\{\lastr(n,e),(n-1)/2\}$ such that
$\binom{r}{2}e+r \ge n-1$, if $h(E_r)\ge 3$ then
\begin{equation*}
  4r^2+2r+e(2en-3e-n)\le \sqrt{3}e(n-1)(2e-1).
\end{equation*}
Moreover, there exists an integer
$r\le \lastr(n,e)$ such that $h(E_r)\le 2$.
\end{Lem}
 
\begin{proof} 
As $\binom{r}{2}e+r \ge n-1$,
it follows from lemma \ref{clm:rmin}, claim \ref{bigrcodim} that
$\#(E_r \setminus E_1)\ge e(n-1)/2$.
On the other hand, we also have
\begin{align*}
 \#(E_r \setminus E_1)=&
 \sum_i i\left(\hat \ell_{E_1}(i-1)-\hat \ell_{E_r}(i-1)\right),\\
 \ell(E_1)-\ell(E_r)=&
 \sum_i \left(\hat \ell_{E_1}(i-1)-\hat \ell_{E_r}(i-1)\right),
\end{align*}
and by \ref{clm:steptypes}, claim \ref{steplength}, 
$\Delta_\ell:=\ell(E_1)-\ell(E_r)= 2\mu_1+3\mu_2+\dots+r\mu_{r-1}.$
Write $\Delta_\ell =\kappa (n-1)+\rho$, with $0\le
\rho<n-1$. As $\hat \ell_{E_1}(i-1)-\hat
\ell_{E_r}(i-1) \le n-1$ for all $i$, an upper bound for
$\#(E_r \setminus E_1)$ is
\begin{equation}
  \label{eq:delta}
  \#(E_r \setminus E_1)\le
\sum_{i=e-\kappa+1}^{e}i(n-1)+(e-\kappa)\rho \le e \Delta_\ell  - 
\frac{\Delta_\ell}{2} \left(\frac{\Delta_\ell}{n-1}-1 \right).
\end{equation}
In particular $e \Delta_\ell  \ge \#(E_r \setminus
E_1)\ge e(n-1)/2$. On the other hand, the function 
$(x/2)(x/(n-1)-1)$ is increasing for $x\ge (n-1)/2$, so 
$$
\frac{\Delta_\ell}{2} \left(\frac{\Delta_\ell}{n-1}-1 \right)\ge
\frac{\#(E_r \setminus E_1)\big/e}{2} \left(\frac{\#(E_r \setminus
    E_1)\big/e}{n-1}-1 \right)
$$
which combined with (\ref{eq:delta}) gives, denoting as before
$M=\mu_1+\mu_2+\dots+\mu_{r_1-1}$,
\begin{equation}
  \label{eq:slices}
\Delta_\ell  \ge  \left(M+\frac{\binom{M+1}{2}}{e}\right)\left(1-
 \frac{1}{2e} +\frac{eM+\binom{M+1}{2}}{2e^2(n-1)}\right).
\end{equation}
On the other hand, claim \ref{steptotal} of theorem
\ref{clm:steptypes} tells us that
\begin{gather*}
  E_r \supset \tau(E_{r-1},(r-1)\mu_{r-1}) \supset \\
\supset \tau(\tau(E_{r-2},(r-2)\mu_{r-2}),(r-1)\mu_{r-1}) =
\tau(E_{r-2},(r-2)\mu_{r-2}+(r-1)\mu_{r-1}) \supset \\
\supset \dots \supset \tau(E_{1},\small\sum
i\mu_i)=\tau(E_{1},\Delta_\ell-M).
\end{gather*}
Therefore 
$$h_{E_r}\le h_{\tau(E_{1},\Delta_\ell-M)}=h_{E_1}(\Delta_\ell-M)
\le e-\frac{\Delta_\ell-M-(n-2)}{n-1},$$
which solving for $\Delta_\ell$ and combining with (\ref{eq:slices}),
gives the bound
\begin{equation}
  \label{eq:hfromM}
 h_{E_r}\le e-\frac{eM+\binom{M+1}{2}}{e(n-1)}\left(1-
 \frac{1}{2e} +\frac{eM+\binom{M+1}{2}}{2e^2(n-1)}\right)+
\frac{M+n-2}{n-1}, 
\end{equation}
where the expression on the right is a decreasing function of $M$.
On the other hand, because $(L-MD,E_r,r)$ is consistent and
$\ell_{E_r}(h(E_r)-1)\ge r$, it follows
 that $e+M>r(h(E_r)-1)$; plugging the resulting bound on $M$ into
 (\ref{eq:hfromM}) and simplifying we get the following
 inequality:
\begin{gather*}
  \left(r^2(h(E_r)-1)^2+r(h(E_r)-1)+e(2en-3e-n)\right)^2\le\\
\le e^2(n-1)^2\left(12e^2+4e+1-8e\frac{e+r+1+(n-1-r)h(E_r)}{n-1}\right).
\end{gather*}
Now, if for some $r\le(n-1)/2$, $h(E_{r})\ge 3$ the
previous inequality gives
\begin{equation}
\label{eq:boundrnm}
  4r^2+2r+e(2en-3e-n)\le \sqrt{3}e(n-1)(2e-1). 
\end{equation}
Finally, recall from the proof of lemma \ref{clm:rmin} that because 
$\lmin(E_r)\ge n-r$,
$r_0=\lfloor (n-e+1)/2 \rfloor \le \lastr(n,e)$.
We shall prove that $h(E_{r_0})\le 2$, by contradiction: suppose
$h(E_{r_0})\ge 3$. 
Since $\binom{r_0}{2}e+r_0 \ge n-1$ and $2r_0 \le n-1$,
$r_0$ must satisfy the inequality (\ref{eq:boundrnm}). Then using
$r_0\ge(n-e)/2$ and the hypothesis $n\ge 4e^2$, we end up with
\begin{equation*}
e\,\left(2\,e -1\right) \,
  \left(  4\,( 3 - {\sqrt{3}} ) \,e^2+1 + {\sqrt{3}} 
    \right)\le 0,
\end{equation*}
which is absurd for positive $e$.
\end{proof}

The following result rephrases Theorem 4.3 in
\cite{Roe06}, and its corollary below follows by iteration, 
in an analogous way to corollary 4.6 in \cite{Roe06}.

\begin{Teo}
\label{clm:below2}
  Let $E$ be a staircase of height two and $s$ a positive integer
  satisfying $\hat \ell_E(0)\ge s+2$ and $\ell_E(1) \ge s$, and 
  let $L$ be a divisor class with $L\cdot D=2s$. Define
  $E_1$ to be the unique staircase of height (at most) two with
 $\ell(E_1)= \ell(E)-s-1$, $ \ell_{E_1}(1)=\ell_E(1)-s$. If
  $\ell_E(1) \ge 2s-1$, define furthermore $E_2$ to be the unique
  staircase of height (at most) two with 
 $\ell(E_2)= \ell(E)-2s-2$, $ \ell_{E_1}(1)=\ell_E(1)-2s-1$.

$(L,E,s)$ is a consistent type.
If $\ell_E(1) \le 2s-2$ then $(L-D,E_1,s+1)$ is a consistent type, and
every linear system of type
$(L-D,E_1,s+1)$ contains as a sublinear system the moving part of a limit of linear systems of type
$(L,E,s)$. If $\ell_E(1) \ge 2s-1$ then $(L-2D,E_2,s+1)$ is a
consistent type, every linear system of type
$(L-2D,E_2,s+1)$ contains as a sublinear system the moving part of a
limit of linear systems of type $(L,E,s)$. 
\end{Teo}

\begin{Cor}
\label{clm:kstepsinone}
    Let $E$ be a staircase of height two and $s$, $k$ positive integers
  satisfying $\hat \ell_E(0)\ge s+2k-2$, $\ell(E)\ge (2s+k)(k-1)$ and
  $\ell_E(1) \ge (2s+k-1)(k-1)$, and  
  let $L$ be a divisor class with $L\cdot D=2s$. Define
  $E_1=E$ and for $k>1$ let $E_k$ be the unique staircase of height
  (at most) two with 
 $\ell(E_k)= \ell(E)-(2s+k)(k-1)$, $\ell_{E_k}(1)=\ell_E(1)-(2s+k-1)(k-1)$. 

$(L-2(k-1)D,E_k,s+k-1)$ is a consistent type, and if there is a
linear system of this type with the expected dimension then there is a
linear system of type $(L,E,s)$ with the expected dimension.
\end{Cor}

The regular linear systems to which everything else is specialized are
those of the following lemma, equivalent to 
lemma 4.4 of \cite{Roe01a}.

\begin{Lem}
\label{clm:bottom}
  Let $E$ be a staircase of height two and $c$ a positive integer
  satisfying $\ell(E)>c$ and $2\ell_E(1)\le c$. Then for every divisor 
class $L$ with $L\cdot D=c$ and every integer $r$ such that the type 
$(L,E,r)$ is consistent, general linear systems of type $(L,E,r)$ are regular.
\end{Lem}

\begin{Lem}
\label{clm:regheight2}
  Let $E$ be a staircase with $h(E)=2$. Let $c$ be a positive integer
  such that
$$ \hat\ell_E(0)+ \frac{c}{2}> 1+
3\sqrt{\ell_E(1)+\left(\frac{c}{2}\right)^2}.$$
Then for every divisor class $L$ with $L\cdot D=c$ and every integer
$r$ such that the type $(L,E,r)$ is consistent, general linear systems
of type $(L,E,r)$ are regular.
\end{Lem}
\begin{proof}
Note that $\ell(E)\ge \hat \ell_E(0)\ge 1-c/2+3c/2>c$, so
if $2\ell_E(1)\le c$ then the result follows from lemma
\ref{clm:bottom}; we assume from now on that \etag{l1gt} $2\ell_E(1)\ge c+1$.
Also, if general linear systems of type $(L,E,r+1)$ are consistent and
regular then general linear systems of type $(L,E,r)$ are consistent
and by semicontinuity regular as well.

For technical reasons we suppose that $c=2s$ is even; let
us see that this is not restrictive. If  $c=2s+1$ is odd then
\eqref{l1gt} gives $\ell_E(1)\ge s+1$ and the hypothesis of the lemma implies
$\hat\ell_E(0)\ge 2s+1$, so for $r\le s+1$, $h^r_E=2$. Thus corollary
\ref{clm:indepf} gives $\len(Z_{p,E,C}\cap D)=2r$ for all $C$ such
that $(C\cdot D)_p=r$. In particular, $(L,E,r)$ is consistent if and
only if $r\le s$. Now we may apply theorem \ref{clm:stcres} or just
specialize to $Z_{p,E,C}$ with $(C\cdot D)_p=s+1$ and we obtain
that it is enough to prove that general linear systems
of type $(L-D,\tau(E,s+1),s+1)$ are regular. But $E'=\tau(E,s+1)$ and
$c'=c+1$ satisfy 
$$ \hat\ell_{E'}(0)+ \frac{c'}{2}> 1+
3\sqrt{\ell_{E'}(1)+\left(\frac{c'}{2}\right)^2}$$
and $(L-D)\cdot D=c'$ so we have reduced to a case with even $c'$.

So we can assume \etag{evl} $c=2s$, and it is not hard to see that also
in this case $(L,E,r)$ is consistent if and
only if $r\le s$. Let $k$ be the integer such that 
$$(k+s)^2 \overset{\etag{kgt}}> \ell_E(1)+s^2\overset{\etag{klt}}\ge (k+s-1)^2$$
(in particular $k\ge 1$). 
The hypothesis of the lemma gives \etag{l0gt} $\hat \ell_E(0)\ge 1+3(k+s-1)-s=2s+3k-2$.

\eqref{klt} can be rewritten as \eqrefprime{klt} $\ell_E(1) \ge (2s+k-1)(k-1)$, which
added to \eqref{l0gt} gives \etag{lgt} $\ell(E)\ge (2s+k+1)k-1$.
Corollary \ref{clm:kstepsinone} has weaker hypotheses than
the inequalities \eqref{l0gt}, \eqref{lgt} and \eqrefprime{klt}, so it applies to 
the present situation. Thus we are reduced to proving that general 
  linear systems of type $(L-2(k-1)D,E_k,s+k-1)$ are regular. Note
  that \eqref{l0gt} and $s+k\ge 2$ imply \etag{l0big}
  $\hat \ell_{E_k}(0)\ge s+k+1$.

We distinguish two cases. If $\ell_E(1) \le (2s+k-1)k-s$ then
$\ell_{E_k}(1)\le s+k-1$, and \eqref{lgt} gives $\ell(E_k)\ge 2s+2k-1$. So lemma \ref{clm:bottom}
proves the needed regularity.

Alternatively, if $\ell_E(1) \ge (2s+k-1)k-s+1$ then \etag{l1big}
  $\ell_{E_k}(1) \ge s+k$. Now theorem 
 \ref{clm:below2} applies (with $s'=s+k-1$) due to \eqref{l0big} and
 \eqref{l1big}, so 
 it is enough to prove that general linear systems of type
 $(L-(2k-1)D,E_k',s+k)$ are regular,  where $E_k'$ has
  height two and $\ell(E_k') =\ell(E_k) -(s+k)$, $\ell_{E_k'}(1)
  =\ell_{E_k}(1)- (s+k-1)$. Adding \eqref{l0big} and \eqref{l1big}
  gives $\ell(E_k') \ge 2k+2s+1$ and \eqref{kgt} implies 
  $\ell_{E_k'}(1)\le k+s-1$, so lemma \ref{clm:bottom}
proves the needed regularity again.
\end{proof}
We are now ready to prove our main theorem.
\begin{proof}[Proof of \ref{clm:smallmult}]
    It is well known, and follows from results of \cite{Roe01},
  that by a semicontinuity argument one can restrict to
  the case that $\D$ is a unibranched diagram of exactly $n$
  free vertices of multiplicity $e$. Moreover, this case is equivalent
  to proving that, taking $\pi:S\rightarrow \P^2$ to be the blow-up of
  the plane at a point, $D$ the exceptional divisor, and $L=\pi^*H-eD$,
  where $H$ is the class of a line, general linear systems on $S$ of type
  $(L,E_1(n,e),1)$ are regular. 

By theorem \ref{clm:steptypes}, every linear system of type
$(L_2,E_2(n,e),2)$ contains as a sublinear system the moving part of a
limit of linear systems of type $(L,E_1(n,e),1)$ and their expected
dimensions agree by claim \ref{stepcodim} of theorem \ref{clm:steptypes} and
\cite[2.14]{Del73}, so it will be enough to show that general linear 
systems of type $(L_2,E_2(n,e),2)$ are regular. Iterating the process,
it is enough to prove that general linear
systems of type $(L_r,E_r(n,e),r)$ are regular, for some $r\le \lastr(n,e)$.

Let $r$ be the minimal integer such that 
$h(E_r)\le 2$. Such an $r$ exists by lemma \ref{clm:goesto2}. Applying
\ref{clm:regheight2}, it will be enough to show that 
\begin{equation}
  \label{eq:goal}
 \hat\ell_{E_r}(0)+ \frac{L_r \cdot D}{2}> 1+
3\sqrt{\ell_{E_r}(1)+\left(\frac{L_r \cdot D}{2}\right)^2}. 
\end{equation}
Now, due to theorem \ref{clm:steptypes}, claim \ref{stepcodim},
$$\ell_{E_r}(0) + \ell_{E_r}(1)+\binom{L_r \cdot D+1}{2}=n\binom{e+1}{2},$$
so (\ref{eq:goal}) is equivalent to
\begin{equation*}
 \hat\ell_{E_r}(0)+ \frac{L_r \cdot D}{2}> 1+
3\sqrt{\frac{n}{2}\binom{e+1}{2}-\frac{\hat\ell_{E_r}(0)}{2}-\frac{L_r
    \cdot D}{4}},
\end{equation*}
and moreover, since $\ell_{E_r}(0) + \ell_{E_r}(1)\le 3(n-1)$ it
follows that $\binom{L_r \cdot D+1}{2}\ge n\binom{e+1}{2}-3(n-1)$ which,
taking into account that $e>2$, implies $L_r \cdot D \ge
e\sqrt{n}-1/2$, so it will be enough to prove
\begin{equation}
  \label{eq:goal2}
 \hat\ell_{E_r}(0)+ \frac{2e\sqrt{n}-1}{4}> 1+
3\sqrt{\frac{n}{2}\binom{e+1}{2}-\frac{\hat\ell_{E_r}(0)}{2}-\frac{2e\sqrt{n}-1}
  {8}}.
\end{equation}
But by claim \ref{stepstairs} of theorem
\ref{clm:steptypes}, $\hat \ell_{E_r}(0)\ge n-r$ 
and the minimality of $r$ together with lemma \ref{clm:goesto2} give
\begin{equation}
\label{ineqlast}
  4(r-1)^2+2(r-1)+e(2en-3e-n)\le \sqrt{3}e(n-1)(2e-1).
\end{equation}
It is now a simple calculus exercise to check that if $e$, $n$, $r$ and
$\hat\ell_{E_r}(0)$ are integers satisfying $e>2$, 
$n\ge 4e^2$, $\hat\ell_{E_r}(0) \ge n-r$ and \eqref{ineqlast}, then
\eqref{eq:goal2} holds.
\end{proof}

\label{sec:hhsquares}

Finally we prove \'Evain's result for a square number of points.

\begin{proof}[proof of \ref{clm:squares}]
It is known after \cite{Eva99} and \cite{BZ03} that the result is true
for $n$ a power of four or nine, and that if it is true for $n_1$ and $n_2$,
then it is true for $n_1 n_2$. So it is not restrictive to assume that
$s$ is odd and $s\ge 5$. On the other hand, due to theorem
\ref{clm:regsmalle}, we may assume that $e>s/2$. 

Let $p_1, \dots, p_{s^2} \in \P^2$ be points in general position,
let $Z$ be the union of these points taken with multiplicity $e$, and
$\I_Z$ the defining ideal sheaf.
By \cite[5.3]{HHF03} it is known that if $a\ge se+(s-3)/2$, then 
$H^1(\P^2,\I_Z(d))=0$ (the map $\rho_{n,e}(d)$ is surjective), and by
\ref{clm:sqodd} we know that  
if $a\le se+(s-5)/2$, then $H^0(\P^2,\I_Z(d))=0$ (the map
$\rho_{n,e}(d)$ is injective) so we are done.
\end{proof}

\bibliographystyle{amsplain}
\bibliography{Biblio}
\end{document}